\documentclass[12pt]{article}

\usepackage{amsmath}
\usepackage{graphicx,psfrag,epsf}
\usepackage{enumerate}
\usepackage{graphicx}
\usepackage{amsmath, amsfonts, amssymb, mathrsfs, rotating, setspace}
\usepackage{textcomp}
\usepackage{enumitem}
\usepackage[utf8]{inputenc}
\usepackage[margin=1in]{geometry}
\usepackage{url} % not crucial - just used below for the URL 
\usepackage{algorithm,algorithmic}
\usepackage{natbib}
\usepackage{xcolor}
\usepackage{bbm}
\usepackage{rotating}
\usepackage{booktabs}
\usepackage{caption}
\usepackage{subcaption}
\usepackage{hyperref}
\usepackage{url}
\usepackage{natbib}
\usepackage{breakurl}
\usepackage{tablefootnote}
\usepackage{booktabs,caption}
\usepackage{tikz}
\usepackage{soul}

\newcommand{\cithree}[3]{_{{#1}\ \!}{#2}_{\ {#3}}}

\usetikzlibrary{positioning}
\usetikzlibrary{arrows}

\usepackage[flushleft]{threeparttable}

\newcommand{\blind}{1}

\begin{document}

\def\spacingset#1{\renewcommand{\baselinestretch}%
{#1}\small\normalsize} \spacingset{1}

%%%%%%%%%%%%%%%%%%%%%%%%%%%%%%%%%%%%%%%%%%%%%%%%%%%%%%%%%%%%%%%%%%%%%%%%%%%%%%

\if1\blind
{
  \title{\bf Nonparametric bounds for causal effects in imperfect randomized experiments}
  \author{Erin E. Gabriel$^1$\thanks{EEG is partially supported by Swedish Research Council grant 2017-01898, AS by Swedish Research Council grant 2016-01267 and MCS by Swedish Research Council grant 2019-00227} \hspace{.2cm}\\
  Arvid Sj\"olander$^1$\\
  Michael C. Sachs$^1$\\
\footnotesize 1:Department of Medical Epidemiology and Biostatistics, Karolinska Institutet, Stockholm, Sweden \\

corresponding author: erin.gabriel@ki.se}
  \maketitle
} \fi 

\if0\blind
{
  \bigskip
  \bigskip
  \bigskip
  \begin{center}
    {\LARGE\bf Nonparametric bounds for causal effects in imperfect randomized experiments}
\end{center}
  \medskip
} \fi

\bigskip

\begin{abstract}
Nonignorable missingness and noncompliance can occur even in well-designed randomized experiments making the intervention effect that the experiment was designed to estimate nonidentifiable. Nonparametric causal bounds provide a way to narrow the range of possible values for a nonidentifiable causal effect with minimal assumptions. We derive novel bounds for the causal risk difference for a binary outcome and intervention in randomized experiments with nonignorable missingness caused by a variety of mechanisms and with or without noncompliance. We illustrate the use of the proposed bounds in our motivating data example of peanut consumption on the development of peanut allergies in infants.
\end{abstract}

\noindent%
\vfill

\newpage
\spacingset{2}
\vspace{-1cm}
\section{Introduction}
The goal of randomized experiments is to estimate the causal effect of an intervention such as a medical treatment, vaccine, or social program. However, when the sample arrived upon at the end of the study is missing outcome information, the causal effect may be nonidentifiable. When there is no missing data, randomization allows for the identification of the the effect of being assigned to the intervention, sometimes called the intent to treat (ITT) effect; this is only equivalent to the intervention effect if subjects comply with their assigned intervention as directed. When this is not the case the intervention effect can also be nonidentifiable, even with no missing data.

There are few papers that focus on bounding nonidentified causal effects in randomized experiments with missing data. A notable exception is \citet{horowitz2000nonparametric} who derive bounds for the risk difference conditional on a measured baseline covariate, making no assumptions about the missingness mechanism. \citet{marden2018implementation} derive bounds for population proportions under nonignorable missing outcome data, but not causal contrasts. Additionally, practitioners almost always calculate an assumption free bound when outcome data are missing in a trial by imputing missing data in the least favourable way for the intervention. Specifically, if the intervention is expected to reduce the probability of the outcome being equal to 1, missing outcomes in the intervention arm would be imputed as 1, and in the control arm as 0, which is recommended as a sensitivity analysis by \citet{ema2011missing}. 
One can form bounds by additionally imputing in the most favourable way possible obtaining what we will call the best/worst case bounds.

Noncompliance is a well known concept in the causal inference literature. \citet{Balke97} developed nonparametric bounds for the causal risk difference when subjects may not comply with the assigned intervention. When noncompliance is compounded by missing outcome data due to study drop-out, loss to follow-up and withdrawal of consent, the standard method of best/worst case imputation does not bound the intervention effect. To our knowledge, bounds for the intervention effect have not yet been derived for settings with both nonignorable missingness and noncompliance. 

Much of the nonparametric causal bounds literature uses the method developed in \citet{balke1994counterfactual} for deriving valid and tight bounds. Valid means that there are no values of the true causal effect outside of the bounds, while tight means that there are no values inside the bounds that the true causal effect can not take on given the available information and assumptions. In order to use this method, the causal effect of interest and the constraints implied by the causal model must be stated as a linear optimization problem. For this reason, much of the literature on nonparametric bounds for causal effects has focused on simple random sampling in observational studies and completely observed data in randomized experiments, which can be easily stated as linear programming problems provided the causal target is linear. \citet{kuroki2010sharp} and \citet{Gabriel2020} are exceptions who derive bounds in settings that are nonlinear. \citet{kuroki2010sharp} derive bounds for the risk ratio under case-control and cohort sampling with and without missing exposure data. \citet{Gabriel2020} derive bounds under more general outcome-dependent observational studies. Although nonignorable missingness can be considered a form of outcome-dependent sampling, \citet{Gabriel2020} and \citet{kuroki2010sharp} do not consider settings with randomized exposure. 

We derive bounds for the causal risk difference of an intervention under a variety of settings with nonignorable missingness of the outcome, with and without noncompliance, which is also subject to missingness, in randomized experiments. We consider three settings with perfect compliance, with differing forms of nonignorable missing data, and five settings that also have noncompliance. We only consider settings where missingness would make observation of compliance impossible, such as in our motivating example, where the intervention (peanut exposure) occurs repeatedly over long-term follow-up up to the time of the outcome measurement. While all three settings we consider under perfect compliance are novel, to our knowledge, three of the five scenarios we consider in the noncompliance settings are equivalent to instrumental variable scenarios considered in \citet{Gabriel2020}. In addition, in settings with noncompliance and nonignorable missingness of the outcome, we provide novel bounds under the assumption of no defiers, which in some settings are tighter than the bounds not assuming no defiers. 

We map each of the scenarios with noncompliance to a scenario with perfect compliance to consider bounds for the ITT or assignment effect, which is then comparable to the best/worst case bounds in those settings, as the best/worst case bounds are for the assignment effect and not the intervention effect. Because of this difference in estimand, the best/worst imputation, which is often considered the most robust or least biased way to report effects in imperfect trials, can actually give much narrower bounds that do not even contain the causal effect of intervention when ignoring noncompliance. For this reason, when compliance is assessed, we recommend using our proposed bounds for the intervention effect in addition to best/worst case imputation for the assignment effect.  

In our motivating example of the regular exposure of infants to peanut products prior to 60 months of life on allergic reactions to peanuts at 60 months there is both observed noncompliance and missingness due to dropout. In the primary publication of this trial, the classic worst case imputation method is used as a sensitivity analysis to nonignorable missingness \citep{du2015randomized}. We demonstrate that although this procedure covers the assignment effect, there is much greater uncertainty in the causal risk difference for the intervention. However, as all bounds exclude a null effect, we strongly confirm the findings of the study that regular exposure of infants to peanut products reduces their risk of peanut allergies later in life. 

The paper is structured as follows, in Section \ref{sec:not} we define notation, provide basic definitions, assumptions, describe the causal models of interest, and review the relevant previously derived bounds. In Section \ref{methods} we describe the methods that we use to derive the novel bounds that we present in Section \ref{newbounds}. In Section \ref{refineb} we qualitatively compare the novel bounds and in Section \ref{num} we carry out a simulation study to assess their performance. In Section \ref{real} we analyse and discuss our motivating example, before providing a summary and discussion of future work and limitations in Section \ref{dis}.

\section{Preliminaries} \label{sec:not}

\subsection{Notation}
Let $X$ be the binary intervention, $Y$ the binary outcome of interest, and $Y(x)$ be the potential (or counterfactual) outcome \citep{rubin1974estimating,pearl2009causality} for a given subject, if the intervention is set to level $x$. Let $O$ be an indicator of having observable outcome and compliance information; $O=1$ for ``observable'' and $O=0$ for ``not observable''. Let $U$ be a set of unobserved variables that will represent common causes or confounders with no restrictions on the distribution of $U$. Thus, the observed data distribution is given by $p\{X,Y|O=1\}$; $p\{\cdot\}$ denotes the probability mass function. As this is a randomized trial and we know all subjects' $X$ values we observe  $p\{O=1|X=x\}$, and therefore the probabilities of interest $p\{Y=y,O=1|X=x\}=p\{Y=y|O=1,X=x\}p\{O=1|X=x\}$ are observable or estimable. 

When compliance is imperfect, the randomization and the actual intervention are not the same. Let $R$ be the assignment of a subject to $X$, which is always randomized with $R=1$ meaning that one was randomized to $X=1$, and $R=0$ to $X=0$. Let $Y(r)$ be the potential (or counterfactual) outcome for a given subject, if the randomization is set to level $r$, and let $X(r)$ be the same for the intervention. In this setting we observe $p\{X,Y,R|O=1\}$, but because we are only considering randomized trials, one will also always know the marginal probabilities of $p\{R=1\}$ and $p\{O=1\}$ and therefore, $p\{O=1|R=r\}$. We can use this to obtain the probabilities of interest in this setting
\begin{eqnarray*}
&&p\{X=x,Y=y,O=1|R=r\}=\\
&&\frac{p\{R=r,X=x,Y=y|O=1\}}{\sum_{x,y}p\{R=r,X=x,Y=y|O=1\}}p\{O=1|R=r\}.
\end{eqnarray*}
Note that in the noncompliance setting one may, in some settings, be able to observe $p\{X\}$, however, we do not consider these situations. We also do not consider settings where $X$ may be missing for more subjects than $Y$; we consider a single missingness mechanism where both $X$ and $Y$ are both observed or both missing. 

Our target parameter of interest is the effect of the intervention as measured by the causal risk difference, $$\theta=p\{Y(X=1)=1\}-p\{Y(X=0)=1\}.$$

Though this is likely the causal estimand of interest, in settings with noncompliance or where compliance is unknown one might also consider what might be referred to as the assignment effect, or the ITT effect, $$\tau=p\{Y(R=1)=1\}-p\{Y(R=0)=1\}.$$ 

For convenience of notation, we define the following probability abbreviations. Let 
\begin{eqnarray*}
p_{y.x1} &=& p\{Y=y|X=x,O=1\}, \\
p_{y1.x} &=& p\{Y=y, O=1|X=x\}, \\
p_{xy1.r} &=& p\{X=x,Y=y,O=1|R=r\},\\
p_{xy.1r} &=& p\{X=x,Y=y|O=1,R=r\},\\
p_{o.x} &=& p\{O=o|X=x\},\\
q_{xy.o} &=& p\{X=x, Y=y|O=o\}.
\end{eqnarray*}

\subsection{Settings} 
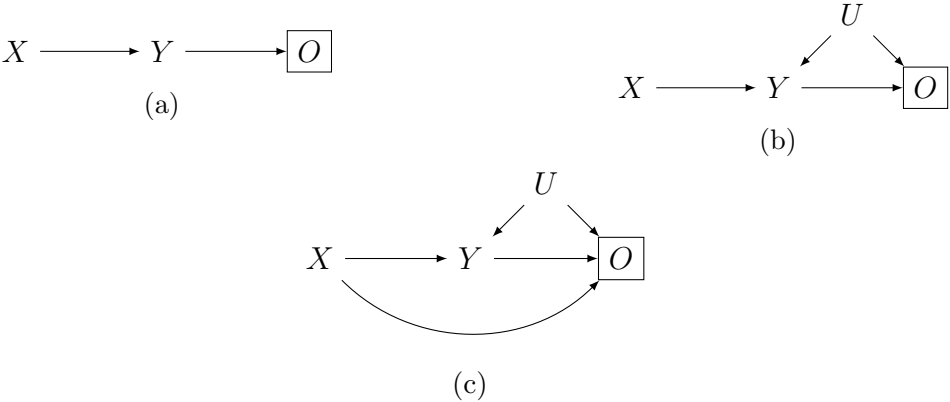
\begin{figure}[ht]

\captionsetup[sub]{width=.9\linewidth}
\centering
\resizebox{\linewidth}{!}{
\begin{subfigure}{0.5\textwidth}
\centering
\begin{tikzpicture}
\node (v) at (0,0) {$X$};
\node (i) at (2,0) {$Y$};
\node[draw, rectangle] (s) at (4,0) {$O$};
\draw[-latex] (v) -- (i);

\draw[-latex] (i) -- (s);
\end{tikzpicture}
\caption{\label{1a}}
\end{subfigure}

\begin{subfigure}{0.5\textwidth}
\centering
\begin{tikzpicture}
\node (uh) at (3,1) {$U$};
\node (v) at (0,0) {$X$};
\node (i) at (2,0) {$Y$};
\node[draw, rectangle] (s) at (4,0) {$O$};
\draw[-latex] (v) -- (i);
\draw[-latex] (i) -- (s);
\draw[-latex] (uh) -- (i);
\draw[-latex] (uh) -- (s);

\end{tikzpicture}
\caption{\label{1b}}
\end{subfigure}
}

\centering

\begin{subfigure}{0.5\textwidth}
\centering
\begin{tikzpicture}
\node (uh) at (3,1) {$U$};
\node (v) at (0,0) {$X$};
\node (i) at (2,0) {$Y$};
\node[draw, rectangle] (s) at (4,0) {$O$};
\draw[-latex] (v) -- (i);
\draw[-latex] (i) -- (s);
\draw[-latex] (uh) -- (i);
\draw[-latex] (uh) -- (s);
\draw[-latex] (v) to[out=-45,in=-135] (s);
\end{tikzpicture}
\caption{\label{1c}}
\end{subfigure}
\caption{Causal diagrams for randomized trials with perfect compliance and nonignorable missing data \label{compliance}}
\end{figure}

\begin{figure}[ht]
\captionsetup[sub]{width=.9\linewidth}
\centering
\resizebox{\linewidth}{!}{
\begin{subfigure}[t]{0.5\textwidth}
\centering
\begin{tikzpicture}
\node (E) at (-1,0.75) {$R$};
\node (uh) at (1,1.75) {$U$};
%\node (u) at (-1,-1) {$V$};
\node (v) at (0,0.75) {$X$};
\node (i) at (2,0.75) {$Y$};
\draw[-latex] (E) -- (v);

\node[draw, rectangle] (s) at (3,0.75) {$O$};
%\draw [dashed,<->] (u) -- (E);
%\draw [-latex] (u) -- (v);

\draw[-latex] (v) -- (i);

\draw[-latex] (uh) -- (i);
\draw[-latex] (i) -- (s);

\draw[-latex] (uh) -- (v);
%\draw[-latex] (uh) -- +(4.75,-0.85) (s);
%\draw[-latex] (uh) -- +(2.79,-0.9) (i);

\end{tikzpicture}
\caption{ \label{2a}}
\end{subfigure}

\begin{subfigure}[t]{0.5\textwidth}
\centering
\begin{tikzpicture}
\node (E) at (-1,0) {$R$};
\node (uh) at (1,1) {$U$};
%\node (u) at (-1,-1) {$V$};
\node (v) at (0,0) {$X$};
\node (i) at (2,0) {$Y$};
\draw[-latex] (E) -- (v);

\node[draw, rectangle] (s) at (3,0) {$O$};
%\draw [dashed,<->] (u) -- (E);
%\draw [-latex] (u) -- (v);

\draw[-latex] (v) -- (i);
\draw[-latex] (i) -- (s);
\draw[-latex] (uh) -- (i);
\draw[-latex] (uh) -- (s);

\draw[-latex] (uh) -- (v);
%\draw[-latex] (uh) -- +(4.75,-0.85) (s);
%\draw[-latex] (uh) -- +(2.79,-0.9) (i);

\end{tikzpicture}
\caption{ \label{2b}}
\end{subfigure}
}

\vspace{4mm}

\captionsetup[sub]{width=.9\linewidth}
\centering
\resizebox{\linewidth}{!}{
\begin{subfigure}[t]{0.5\textwidth}
\centering
\begin{tikzpicture}
\node (E) at (-1,0) {$R$};
\node (uh) at (1,1) {$U$};
%\node (u) at (-1,-1) {$V$};
\node (v) at (0,0) {$X$};
\node (i) at (2,0) {$Y$};
\draw[-latex] (E) -- (v);

\node[draw, rectangle] (s) at (3,0) {$O$};
%\draw [dashed,<->] (u) -- (E);
%\draw [-latex] (u) -- (v);

\draw[-latex] (v) -- (i);
\draw[-latex] (uh) -- (s);
\draw[-latex] (uh) -- (i);
\draw[-latex] (i) -- (s);
\draw[-latex] (v) to[out=-45,in=-135] (s);
\draw[-latex] (uh) -- (v);
%\draw[-latex] (uh) -- +(4.75,-0.85) (s);
%\draw[-latex] (uh) -- +(2.79,-0.9) (i);

\end{tikzpicture}
\caption{ \label{2c}}
\end{subfigure}

\begin{subfigure}[t]{0.5\textwidth}
\centering
\begin{tikzpicture}

\node (E) at (-1,0) {$R$};
\node (uh) at (1,1) {$U$};
%\node (u) at (-1,-1) {$V$};
\node (v) at (0,0) {$X$};
\node (i) at (2,0) {$Y$};
\draw[-latex] (E) -- (v);

\node[draw, rectangle] (s) at (3,0) {$O$};
%\draw [dashed,<->] (u) -- (E);
%\draw [-latex] (u) -- (v);

\draw[-latex] (v) -- (i);

\draw[-latex] (uh) -- (i);
\draw[-latex] (uh) -- (s);
\draw[-latex] (E) to[out=-45,in=-135] (s);
\draw[-latex] (uh) -- (v);
\draw[-latex] (i) -- (s);
%\draw[-latex] (uh) -- +(4.75,-0.85) (s);
%\draw[-latex] (uh) -- +(2.79,-0.9) (i);

\end{tikzpicture}
\caption{\label{2d}}
\end{subfigure}
}

\begin{subfigure}[t]{0.5\textwidth}
\centering
\begin{tikzpicture}
\node (E) at (-1,0) {$R$};
\node (uh) at (1,1) {$U$};
%\node (u) at (-1,-1) {$V$};
\node (v) at (0,0) {$X$};
\node (i) at (2,0) {$Y$};
\draw[-latex] (E) -- (v);

\node[draw, rectangle] (s) at (3,0) {$O$};
%\draw [dashed,<->] (u) -- (E);
%\draw [-latex] (u) -- (v);

\draw[-latex] (v) -- (i);
\draw[-latex] (uh) -- (s);
\draw[-latex] (uh) -- (i);
\draw[-latex] (i) -- (s);
\draw[-latex] (v) to[out=-45,in=-135] (s);
\draw[-latex] (uh) -- (v);
\draw[-latex] (E) to[out=-45,in=-135] (s);
%\draw[-latex] (uh) -- +(4.75,-0.85) (s);
%\draw[-latex] (uh) -- +(2.79,-0.9) (i);

\end{tikzpicture}
\caption{ \label{2e}}
\end{subfigure}

\caption{Causal diagrams for randomized trials with noncompliance and nonignorable missing data \label{noncomplicance}}
\end{figure}
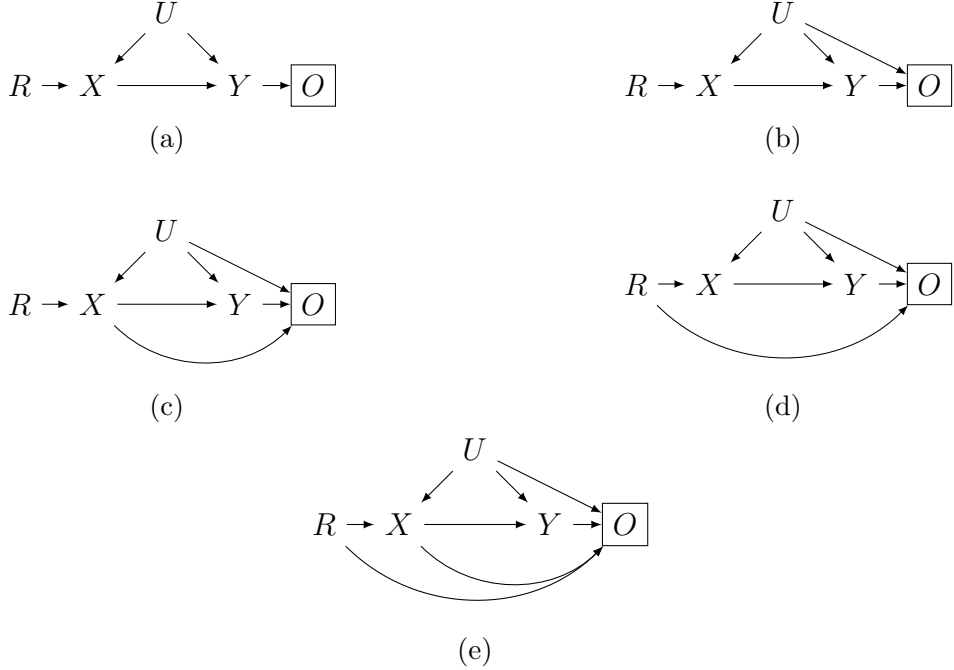

The causal diagrams \citep{pearl2009causality} in Figures \ref{1a} - \ref{2e} represent possible scenarios in a randomized experiment. Figures \ref{1a} - \ref{1c} could be described as randomized experiments with perfect compliance but nonignorable missingness in the outcome. The nonignorable missingness mechanisms we consider are of three types: missingness that is only causally related to the outcome of interest (Figure \ref{1a}), missingness that is associated with the outcome of interest because of an unmeasured common cause of the missingness and the outcome, in addition to being causally related to the outcome (Figure \ref{1b}), and missingness that is additionally causally related to the intervention (Figure \ref{1c}).

Real life settings that fit all the perfect compliance scenarios are single time-point intervention trials where the intervention is administered at the time of randomization. Some examples are a one dose vaccine, a surgical intervention or a single dose intravenous treatment, where subjects may have previous been screened for entry into the study but are not randomized and therefore not actually enrolled until just before the intervention is performed. Although this type of randomization procedure reduces or even eliminates compliance issues, unless the endpoint is immediate, such trials can still suffer from nonignorable missingness in the outcome. In contrast, any time an intervention requires active participation from the subjects under study, compliance as well as missingness can be issues.

The actual intervention $X$ may differ from the randomized assignment $R$, and therefore $X$ and $Y$ are confounded in all settings of Figure 2. Noncompliance alone alone can cause Figure \ref{2a} has noncompliance in addition to nonignorable missingness due to a causal effect of the outcome on the missingness without confounding. Figure \ref{2b} is the same as Figure \ref{2a}, with noncompliance in addition to nonignorable missingness due to a causal effect of the outcome on the missingness, but with additional confounding. Figure \ref{2c} through Figure \ref{2e} depict various causal relationships between the missingness and the outcome, the randomization and the true intervention $X$, but all have nonignorable missingness due to unmeasured common causes of the missingness and the outcome as well as potential causal effects of the outcome, interventions and randomization, all under noncompliance. 

Real life trials that fit Figure 2 include any take-at-home medications, diet or physical activity interventions. When such a trial uses an intervention that is available to all participants, and is not blinded to participants, any type of noncompliance is possible. For example, in a randomized trial of diet and exercise it might be the case that being told not to exercise or diet may induce some participants to exercise, while telling those same subjects to exercise might overwhelm them or make them defensive, thus inducing them to not perform the randomized intervention regardless of their randomization. For this reason, bounds not considering any further assumptions about the type of compliance may be needed in many experimental settings. 

In any of the settings with noncompliance it may be of interest to further consider if it is possible that subjects randomized to a particular intervention would defy it. This assumption can be stated in terms of the counterfactuals as
\begin{eqnarray}
X(r = 1) \geq X(r = 0). \label{nodef}
\end{eqnarray}
\citet{angrist1996identification} and others have referred to this assumption as monotonicity, but we will use the term \emph{no defiers} for clarity. The no defiers assumption is justified in settings with experimental intervention only available to those randomized to it. Placebo subjects will not have access to the intervention and therefore $X(r = 0) = 0$. This setting implies no defiers, but this is not required for no defiers to be a plausible assumption.

Instead, our real data example offers a less restrictive setting where no defiers is plausible, but some randomized to no intervention are still observed to take some form of the intervention. Our real data example is a trial of peanut exposure for infants where children are randomized either to an intervention of consuming peanut products or to avoid all exposure to peanuts in an unblinded manner. Some parents elected to feed their children peanut products in the avoidance arm and some parents elected to avoid peanuts in the intervention arm. Provided the proportion receiving the intervention of peanut products is higher in the arm randomized to the intervention than in those randomized to no intervention, there are no observable ways to rule out no defiers. It is also hard to imagine a rationale that would compel these parents to do the opposite had they been randomized differently, although it is possible that we simply do not observe enough defiers to detect this pattern. We therefore consider bounds in all settings with noncompliance both with and without the no defiers assumption.

\subsection{Previous bounds}
\citet{robins1989analysis} derived bounds in the setting with noncompliance without missing data, Figure \ref{2a} without $O$. However, \citet{Balke97} showed that those bounds are not tight by deriving new bounds using the linear programming method developed in \citet{balke1994counterfactual}. \citet{Gabriel2020} derived bounds that apply in the settings of Figures \ref{2a}, \ref{2b} and \ref{2c}, when no assumptions are made about defiers, in terms of observational studies with instrumental variables, rather than randomization. The bounds for Figures \ref{2a} and \ref{2b}, as given in \citet{Gabriel2020}, without assuming no defiers, are reproduced in the supplementary material. \citet{kuroki2010sharp} derived bounds in terms of probabilities conditional on $Y$ that are applicable in observational settings with nonrandom sampling and potentially missing exposure information, not randomized settings. 

The worst/best case bounds that are often used in practice can be written in terms of the true probabilities as: \begin{eqnarray}
\label{Beq:0}
p_{1.11}p_{1.1}-p_{1.01}p_{1.0}-p_{0.0}\leq \theta \leq 
p_{1.11}p_{1.1}-p_{1.01}p_{1.0}+p_{0.1}. 
\end{eqnarray}
Replacing $x$ in $p_{y.x1}$ and $p_{o.x}$ with $r$, and ignoring $x$, in the case of noncompliance gives the theoretical best/worst case bounds for the assignment effect $\tau$. We will compare to this theoretical version of the best/worst case bounds in what follows any time we are using the true rather than the estimated probabilities. 

\citet{horowitz2000nonparametric}, as mentioned in the introduction, derived bounds for risk difference conditional on a baseline covariate in randomized settings with missing data, making no assumptions about the missingness mechanism. It can easily be shown that in the special case where there is no baseline covariate that the bounds given in their corollary 1 of Theorem 1 simplify to the best/worst bounds given in \eqref{Beq:0}.  

\section{Methods} \label{methods}
\citet{Gabriel2020} modified the method of \citet{balke1994counterfactual} to apply to a partially observed setting, providing bounds in the settings of Figures \ref{2b} and \ref{2c}, which they referred to as confounded outcome-dependent and confounded exposure- and outcome- dependent settings. However they considered these under the conceptual framework of an instrumental variable and a observational study with unmeasured confounding. We will use a similar approach to derive bounds for Figures \ref{1b}-\ref{1c} and \ref{2d}-\ref{2e}. \citet{Gabriel2020} use a different approach to account for the nonlinear constraint implied by the unconfounded sampling, i.e., the lack of an arrow from $U$ to $O$ in Figure \ref{2a}, and a setting similar to Figure \ref{1a} but with unmeasured confounding between $X$ and $Y$. We will follow a similar approach to derive bounds for Figure \ref{1a} and Figure \ref{2a} assuming no defiers.

\subsection{Linear programming}

In order to use this algorithm to derive bounds that are valid and tight, one must derive linear constraints relating observed probabilities to counterfactual probabilities that are necessary and sufficient for the observed distribution to be in the causal model. In the supplementary material, we use the response function variable representation of the causal model to relate each of the observed probabilities $p\{X = x, Y = y, O = 1\}, x, y \in \{0, 1\}$ for Figure 1, and $p\{X = x, Y = y, R = r, O = 1\}, x, y, r \in \{0, 1\}$ for Figure 2, to counterfactual quantities. Under the settings Figure \ref{1b} and \ref{1c}, the equations are linear, and furthermore, they can be factorized so that all of the linear equations can be written in terms of probabilities of the form $p\{Y = y, O = 1 | X = x\}$. Under the settings Figures \ref{2b} - \ref{2e}, the equations are linear and can be factorized into probabilities of the form $p\{X = x, Y = y, O = 1 | R = r\}$. 

We also show that the target quantity $\theta$ is linear in counterfactual quantities. Treating $\theta$ as the objective function and optimizing it subject to the linear constraints in terms of the observed probabilities is a linear programming problem. Solutions to this problem can be found symbolically by applying Balke's implementation of a vertex enumeration algorithm \citep{balke1994counterfactual, mattheiss1973algorithm}. This gives the bounds on the causal effect of interest as the minimum (maximum) of a list of terms involving only observable probabilities, each of which corresponds to a vertex. This demonstrates that for these problems in Figures \ref{1b}, \ref{1c}, \ref{2b}-\ref{2e}, valid and tight bounds on $\theta$ can be derived symbolically in terms of $p_{y1.x}$ and $p_{xy1.r}$ according to this algorithm.

\subsection{Expansion}
In the settings of Figure \ref{1a} and \ref{2a}, the lack of an arrow from $U$ to $O$ implies that the constraints are nonlinear. We will therefore use a different approach that yields valid but non necessarily tight bounds in these settings. For Figure \ref{1a}, we start with the point identified estimator of $\theta$ in the model with no missingness ($O$ absent from the model). For Figure \ref{2a} we start with the valid and tight bounds under the setting where $O$ is absent and data are fully observed that can be derived with the linear programming method. We then partition those quantities into observed and unobserved parts, by conditioning on the unconfounded variable $O$. Finally, we use the nonlinear constraint to bound the unobserved part, thus producing bounds for the target parameter. Detailed derivations using this approach for Figure \ref{1a} are in the supplementary material, while those for Figure \ref{2a} are in \citet{Gabriel2020}.

\section{Novel bounds} \label{newbounds}
\subsection{Figure 1 bounds for $\theta$}
\noindent \textbf{Result 1}:\\
The bounds given in \eqref{Beq:1} and \eqref{Beq:2} are valid for $\theta$ in the setting of Figure \ref{1a} provided that $A(y,0)=q_{1y.o}/q_{0y.o}$ is not undefined for any value $y$ and $p\{X=0|O=0\}>0$ and $p\{X=1|O=0\}>0$. 

\begin{eqnarray}
\label{Beq:1}
&&\theta\geq 1-\left(p_{0.11}p_{1.1}+p_{1.01}p_{1.0}\right)-\max\left\{\frac{1}{1+A(1,1)},\frac{A(0,1)}{1+A(0,1)}\right\}\times \nonumber \\ 
&&\max\left\{\frac{p_{0.1}}{p\{X=1|O=0\}},\frac{p_{0.0}}{p\{X=0|O=0\}}\right\} 
\end{eqnarray}
and
\begin{eqnarray}
\label{Beq:2}
&&\theta\leq 1-\left(p_{0.11}p_{1.1}+p_{1.01}p_{1.0}\right)-\min\left\{\frac{1}{1+A(1,1)},\frac{A(0,1)}{1+A(0,1)}\right\}\times \nonumber \\ 
&&\min\left\{\frac{p_{0.1}}{p\{X=1|O=0\}},\frac{p_{0.0}}{p\{X=0|O=0\}}\right\}. 
\end{eqnarray}

We give detailed derivations of Result 1 in the supplementary material. It is of note that we show that these bounds are valid, but do not claim that they are tight, in the setting of Figure \ref{1a}. In fact, the bounds are not tight and can be made tighter as discussed in Section \ref{refineb}. 

All bounds that follow, other than those in Result 5, use the modification to the linear programming method of \citet{Balke95} that was first introduced in \citet{Gabriel2020}, for partial observation of the joint probabilities of the data. 

\noindent \textbf{Result 2}:\\
The bounds for $\theta$ given in \eqref{Beq:3} and \eqref{Beq:4} are valid and tight in the settings of Figure \ref{1b}.

\begin{eqnarray}
\label{Beq:3}
\theta\geq \max\left\{\begin{array}{l}
p_{01.0} + p_{11.1} - 1\\
- p_{11.0} + 2 p_{11.1} - 1\\
2 p_{01.0} - p_{01.1} - 1\\
\end{array}\right\},
\end{eqnarray}
and
\begin{eqnarray}
\label{Beq:4}
\theta\leq \min\left\{\begin{array}{l}
- p_{01.1} - p_{11.0} + 1\\
- 2 p_{11.0} + p_{11.1} + 1\\
p_{01.0} - 2 p_{01.1} + 1\\
\end{array}\right\}.
\end{eqnarray}

\noindent \textbf{Result 3}:\\
The bounds for $\theta$ given in \eqref{Beq:5} are valid and tight in the settings of Figures \ref{1c}.

\begin{eqnarray}
\label{Beq:5}
p_{01.0} + p_{11.1} - 1 \leq \theta\leq
- p_{01.1} - p_{11.0} + 1.
\end{eqnarray}

\subsection{Figure 2 bounds for $\theta$}

\noindent \textbf{Result 4}:\\
The bounds for $\theta$ given in \eqref{Beq:6} and \eqref{Beq:7} are valid and tight in the settings of Figure \ref{2c}-\ref{2e}.

\begin{eqnarray}
\label{Beq:6}
\theta\geq \max\left\{\begin{array}{l}
p_{001.1} + p_{111.1} - 1\\
  p_{001.0} + p_{111.1} - 1\\
  p_{001.1} + p_{111.0} - 1\\
  p_{001.0} + p_{111.0} - 1\\
  2 p_{001.1} + p_{011.0} + p_{111.0} + p_{111.1} - 2\\
  2 p_{001.0} + p_{011.1} + p_{111.0} + p_{111.1} - 2\\
  p_{001.0} + p_{001.1} + p_{101.0} + 2 p_{111.1} - 2\\
  p_{001.0} + p_{001.1} + p_{101.1} + 2 p_{111.0} - 2
\end{array}\right\},
\end{eqnarray}
and
\begin{eqnarray}\theta\leq \min\left\{\begin{array}{l}
\label{Beq:7}
  - p_{101.0} - p_{011.0} + 1\\
  - p_{101.0} - p_{011.1} + 1\\
  - p_{101.1} - p_{011.0} + 1\\
  - p_{101.1} - p_{011.1} + 1\\
  - p_{001.0} - p_{101.0} - p_{101.1} - 2 p_{011.1} + 2\\
  - p_{001.1} - p_{101.0} - p_{101.1} - 2 p_{011.0} + 2\\
  - 2 p_{101.0} - p_{011.0} - p_{011.1} - p_{111.1} + 2\\
  - 2 p_{101.1} - p_{011.0} - p_{011.1} - p_{111.0} + 2
\end{array}\right\}.
\end{eqnarray}

\subsection{Figure 2 Bounds for $\theta$ under the no defiers assumption}

\noindent \textbf{Result 5}:\\
Under the no defiers assumption, the bounds for $\theta$ given in \eqref{Beq:8} are valid in the settings of Figure \ref{2a}.

\begin{eqnarray}
\label{Beq:8}
p_{11.11}p\{O=1|R=1\}+p_{00.10}p\{O=1|R=0\}-1 \leq \theta \nonumber \\\leq 1-p_{10.11}p\{O=1|R=1\}-p_{01.10}p\{O=1|R=0\}.
\end{eqnarray}

We derive these bounds by starting with the single term bound given in \citet{Balke95} for the setting of Figure \ref{2a} without missing data under the no defiers assumption, then use the same expansion procedure as described above to arrive at the bounds in \eqref{Beq:8}. These are the first term of the lower/upper bounds for Figure \ref{2a} not assuming no defiers, which are given the supplementary materials. 

\noindent \textbf{Result 6}:\\
Under the no defiers assumption the bounds for $\theta$ given in \eqref{Beq:9} and \eqref{Beq:10} are valid and tight in the settings of Figure \ref{2b}.

\begin{eqnarray}
\label{Beq:9}
\theta\geq \max\left\{\begin{array}{l}
p_{001.0} + p_{111.1} - 1 \\
p_{001.1} - p_{011.0} + p_{011.1} - p_{111.0} + 2 p_{111.1} - 1 \\
2 p_{001.0} - p_{001.1} + p_{101.0} - p_{101.1} + p_{111.0} - 1
\end{array}\right\},
\end{eqnarray}
and
\begin{eqnarray}\theta\leq \min\left\{\begin{array}{l}
\label{Beq:10}
- p_{101.1} - p_{011.0} + 1 \\
- p_{101.0} - 2 p_{011.0} + p_{011.1} - p_{111.0} + p_{111.1} + 1 \\
p_{001.0} - p_{001.1} + p_{101.0} - 2 p_{101.1} - p_{011.1} + 1
\end{array}\right\}.
\end{eqnarray}

\noindent \textbf{Result 7}:\\
Under the no defiers assumption the bounds for $\theta$ given in \eqref{Beq:11} and \eqref{Beq:12} are valid and tight in the settings of Figure \ref{2c}-\ref{2e}. In the setting of Figure \ref{2c} the first term in both the lower and the upper are the only active terms in the bounds i.e. the first terms will always be the max/min of the respective lower/upper set of terms.

\begin{eqnarray}
\label{Beq:11}
\theta\geq \max\left\{\begin{array}{l}
p_{001.0} + p_{111.0} - 1 \\
p_{001.1} + p_{111.0} - 1 \\
p_{001.1} + p_{111.1} - 1 \\
p_{001.0} + p_{111.1} - 1
\end{array}\right\},
\end{eqnarray}
and
\begin{eqnarray}\theta\leq \min\left\{\begin{array}{l}
\label{Beq:12}
1 - p_{101.1} - p_{011.0} \\
1 - p_{101.0} - p_{011.1}  \\
1 - p_{101.0} - p_{011.0}  \\
1 - p_{101.1} - p_{011.1} 
\end{array}\right\}.
\end{eqnarray}

\subsection{Figure 2 bounds for $\tau$}
Under perfect compliance, $R=X$, and therefore all Figure 1 bounds are for both $\theta$ and $\tau$. This is not the case with noncompliance. As the ITT or assignment results are often used in randomized clinical trials regardless of noncompliance or missingness issues, we map the bounds for Figure 1 for $\theta$ to the assignment effects bounds for $\tau$ in Figure 2. 
\\

\noindent \textbf{Result 8}:\\
The bounds for $\theta$ given in \eqref{Beq:1} and \eqref{Beq:2} for Figure \ref{1a} are valid for $\tau$, replacing $X$ with $R$, in the setting of Figure \ref{2a}.

\noindent \textbf{Result 9}:\\
The bounds for $\theta$ given in \eqref{Beq:3} and \eqref{Beq:4} for Figure \ref{1b} are valid and tight for $\tau$, replacing $X$ with $R$, under Figure \ref{2b}, and the bounds for $\theta$ given in \eqref{Beq:5} for Figure \ref{1c} are valid and tight for $\tau$, replacing $X$ with $R$, under Figure \ref{2c}-\ref{2e}.

\subsection{Estimation of bounds}
Up to this point the bounds have only been discussed based on true probabilities. However, all proposed bounds are functions of probabilities that can be estimated by their sample proportions to produce estimated bounds. To account for the statistical uncertainty in the estimates due to sampling we suggest the nonparametric bootstrap \citep{efron1979bootstrap}, which we illustrate the use of in both the simulations and the real data example. 

\section{Refinement and comparison of bounds} \label{refineb}
The bounds derived in Result 1 are valid but not tight. To tighten these bounds, note that any bounds that additionally allow for confounding of either the $X-Y$ or $Y-O$ relationships are valid under Figure \ref{1a} as well, and it is possible that they are sometimes narrower. This is shown in Figure S1 of the supplementary material, where the width of the bounds of Result 1 is compared to the width of bounds allowing $X-Y$ confounding (left panel), and allowing $Y-O$ confounding (right panel). The bounds allowing $Y-O$ confounding are in Result 3, and those allowing $X-Y$ confounding were derived in \citet{Gabriel2020}, and are reproduced in the supplementary material in equations 1 and 2. 

This suggests a simple way to improve the bounds in Result 1, namely to replace them with the alternative bounds whenever they are tighter. Formally, let $l_a$ and $u_a$ denote the lower and upper bounds in Result 1, $l_c$ and $u_c$ be the lower and upper bounds in Result 3, and let $l_f$ and $u_f$ be the lower and upper bounds given in equations 1 and 2 of the supplementary materials. We thus define new bounds for $\theta$ under the causal diagram in Figure \ref{1a} that will be used instead of the bounds in Result 1 for the remainder of the paper: 
\begin{eqnarray}
\label{Beq:refa}
\max(l_a, l_c, l_f) \leq \theta \leq \min(u_a, u_c, u_f).
\end{eqnarray}

This refinement also clearly holds for the bounds for $\tau$ in Figure \ref{2a}, and we will therefore use the bounds in \eqref{Beq:refa} for the $\tau$ assignment effect bounds in Figure \ref{2a}, replacing $X$ with $R$, for the remainder of the paper. 

In addition to the refinement in Figure \ref{1a}, it is easily shown that the bounds in \eqref{Beq:5} are equivalent to the best/worst case bounds in \eqref{Beq:0}. Therefore, whenever one uses the worst/best bounds in settings with perfect compliance one is, in effect, allowing for both confounding of the outcome and the missingness, and the missingness to be influenced by both the outcome and the intervention, as in Figure \ref{1c}. Although we make numerical comparisons in the simulations, the bounds in \eqref{Beq:5} will never differ from the bounds in \eqref{Beq:0}.

\citet{Balke97}, and  \citep{robins1989analysis}, found that assuming ``monotonicity'', which they equate to the no defiers assumption as we present it, in the noncompliance setting  results in a set of bounds that are a subset of the bounds not making the no defiers assumption. As pointed out in \citet{Balke95} because of the structure of the bounds, taking the maximum of the lower bounds terms and the minimum the upper, if the bounds derived under no defiers are valid, tight and a subset of the valid bounds not assuming no defiers, there is nothing gained by using the bounds assuming no defiers, as the only active terms in the bounds derived without assuming no defiers must be those terms given in the no defiers bounds, when there are in fact no defiers. Otherwise, the bounds assuming no defiers would either be invalid or not tight.  

We also find that the tight and valid bounds derived via the linear programming method in the settings of Figure \ref{2c}-\ref{2e} are a subset of the bounds derived without making the no defiers assumption. As they are both tight and valid, this implies that there is again nothing gained by assuming no defiers, as having no defiers will automatically make those terms displayed in \eqref{Beq:11} and \eqref{Beq:12} the only active terms in \eqref{Beq:6} and \eqref{Beq:7}; or, under Figure \ref{2c}, the single terms from the set of four. We also demonstrate this via simulation. This does not hold however, under Figure \ref{2b}, as the terms given in \eqref{Beq:9} and \eqref{Beq:10} are not a subset of those in \eqref{Beq:3} and \eqref{Beq:4}, and are occasionally tighter under no defiers, a fact we demonstrate via simulation. 

The bounds given in \eqref{Beq:6} and \eqref{Beq:7}, which are valid and tight in Figure \ref{2c}-\ref{2e}, become the bounds in \eqref{Beq:5} when $R=X$. A similar equivalence was observed in the noncompliance setting with no missing data in \citep{Balke97}. This is similarly true in the setting of Figure \ref{2b}, although we do not reproduce the bounds not assuming no defiers here. Considering the bounds assuming no defiers given in \eqref{Beq:9} and \eqref{Beq:10}, it is easy to see that if $R=X$, then these bounds become those given for Figure \ref{1b}, \eqref{Beq:3}, \eqref{Beq:4}. In the case of Figure \ref{2a}, the bounds are valid but not tight, thus the connection between the bounds we give for Figure \ref{1a} and the bounds for Figure \ref{2a} under no defiers is not as clear. 

\section{Simulations} \label{num}
We carried out simulation studies in order to compare the width of the true bounds across the different causal diagrams, assess the impact of the amount of missingness on the width of the true bounds, and also to assess the performance of estimated bounds based on samples. For the settings in Figure 1, we generate probability distributions $p\{U,X,Y,O\}$ under the model
\begin{eqnarray}
\label{eq:simmod}
\left .
\begin{array}{rcl}
&&U \sim  \mbox{N}(0,1) \\
%&&p\{R=1\}  \sim  \mbox{Unif}(0.2,0.8) \\
&&p\{X=1\}  =  \mbox{expit}(\alpha_1)\\
&&p\{Y=1|U,X\}  = \mbox{expit}(\beta_1+\delta_1\beta_2U+\beta_3X)\\
&&p\{O=1|U,Y, X\}  =  \mbox{expit}(\gamma_1+\delta_2\gamma_2U +\gamma_3Y + \delta_3\gamma_4X)\\
&&(\alpha_1,\beta_1,\beta_2,\beta_3,\gamma_1,\gamma_2, \gamma_3, \gamma_4)  \sim  N(0,4)\\
&&(\delta_1, \delta_2, \delta_3)   \in  \{0, 1\}
\end{array}
\right \}
\end{eqnarray}
where $\mbox{expit}(x) = e^x / (1 + e^x)$ and where $e$ is Euler's number. The constants $\delta_1, \delta_2, \delta_3$ determine under which of the settings in Figure \ref{1a} - \ref{1c} the distributions are generated: Figure \ref{1a} with $\delta_1 = \delta_2 = \delta_3 = 0$, Figure \ref{1b} with $\delta_1 = \delta_2 = 1, \delta_3 = 0$, and Figure \ref{1c} with $\delta_1 = \delta_2 = \delta_3 = 1$.

For the settings with noncompliance in Figure 2, we generate probability distributions $p\{U,R,X,Y,O\}$ by modifying the model \eqref{eq:simmod} by
\begin{eqnarray}
\label{eq:simmod2}
\left .
\begin{array}{rcl}
&&p\{R=1\}  \sim  \mbox{Unif}(0.2,0.8) \\
&&p\{X=1 |U,R\}  =  \mbox{expit}(\alpha_1 + \alpha_2U + \alpha_3R)\\
%&&p\{Y=1|U,X\}  = \mbox{expit}(\beta_1+\delta_1\beta_2U+\beta_3X)\\
&&p\{O=1|U,R,X,Y\}  =  \mbox{expit}(\gamma_1+\varepsilon_1\gamma_2U +\gamma_3Y + \varepsilon_2\gamma_4X + \varepsilon_3\gamma_5R)\\
&&(\alpha_2,\alpha_3,\gamma_1,\gamma_5)  \sim  N(0,4)\\
&&(\varepsilon_1, \varepsilon_2, \varepsilon_3)   \in  \{0, 1\}
\end{array}
\right \}
\end{eqnarray}
As above, the constants $\varepsilon_1, \varepsilon_2, \varepsilon_3$ determine which of the 5 settings in Figure 2 are satisfied.

\begin{figure}[ht]
\centering
\includegraphics[width = .95\textwidth]{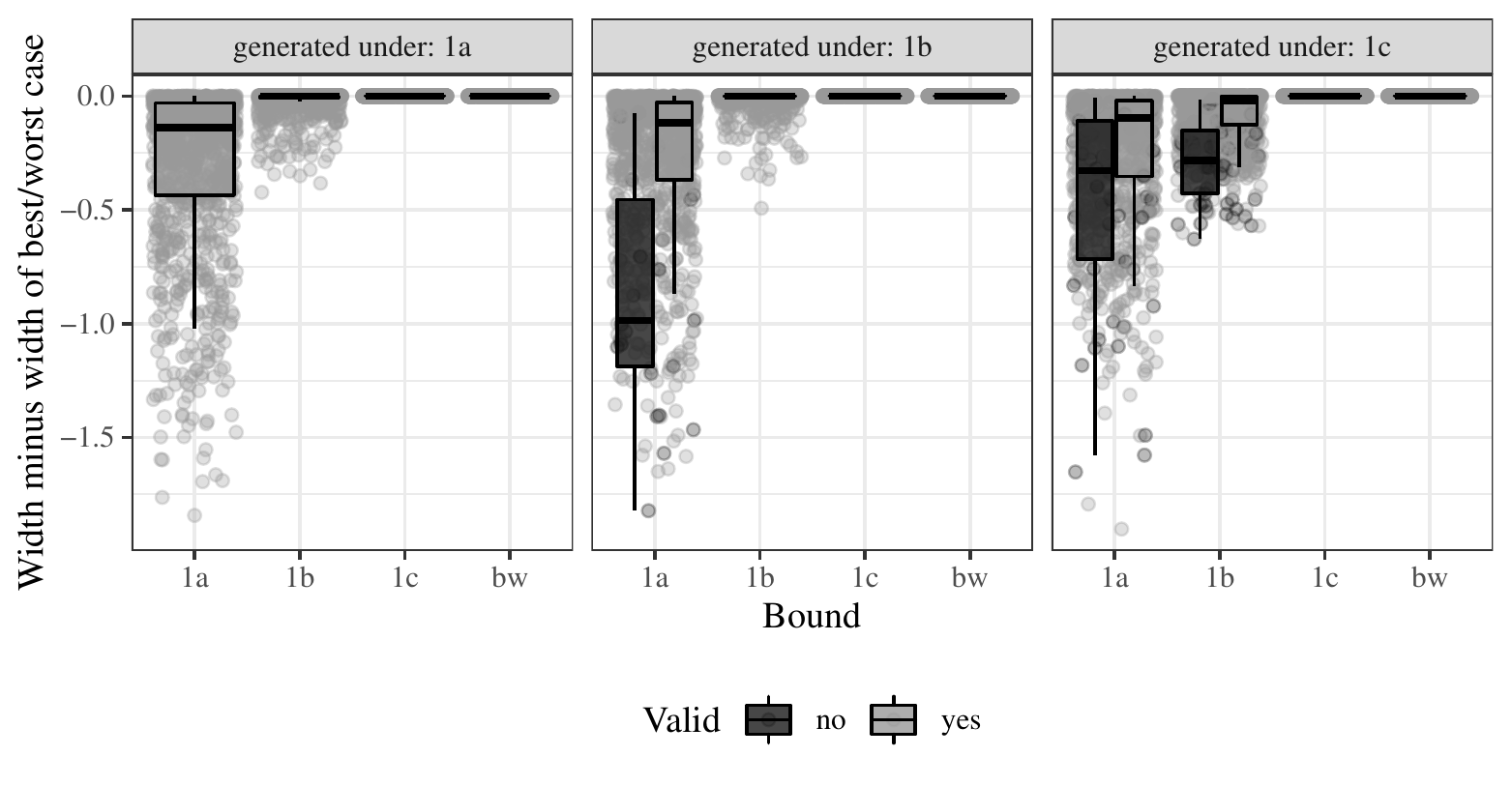}
\caption{\label{relwidth1} Comparison of the width of the true bounds for datasets generated under the DAGs in Figure 1 for distributions that are generated under Figure \ref{1a} (left panel), \ref{1b} (middle panel), and \ref{1c} (right panel). The y-axis shows the width of the bounds for each setting minus the width of the best/worst case bounds (denoted bw in the Figure). The light grey dots and boxes indicate cases where the bounds are valid (i.e., the true value is within the bounds), and the dark grey bounds that are invalid. }
\end{figure}

We first generate 1000 distributions for each setting from the models in \eqref{eq:simmod} and \eqref{eq:simmod2}. Then we compute the bounds under each setting and the best/worst case bounds using the true probabilities generated by the random coefficients. The relative widths of the bounds compared to the best/worst procedure for distributions generated under settings \ref{1a} - \ref{1c} are shown in Figure \ref{relwidth1}. The bounds computed under 1a and 1b are always equal or narrower than the best/worst procedure, however when the distribution does not satisfy setting 1a, the 1a bounds occasionally do not cover the true $\theta$, indicated by darker dots and boxes, and when the distribution does not satisfy 1b, the 1a and 1b bounds occasionally do not cover the truth. The bounds computed under settings 1c are numerically identical to the best/worst procedure, as expected. 

\begin{figure}[ht]
\centering
\includegraphics[width = .95\textwidth]{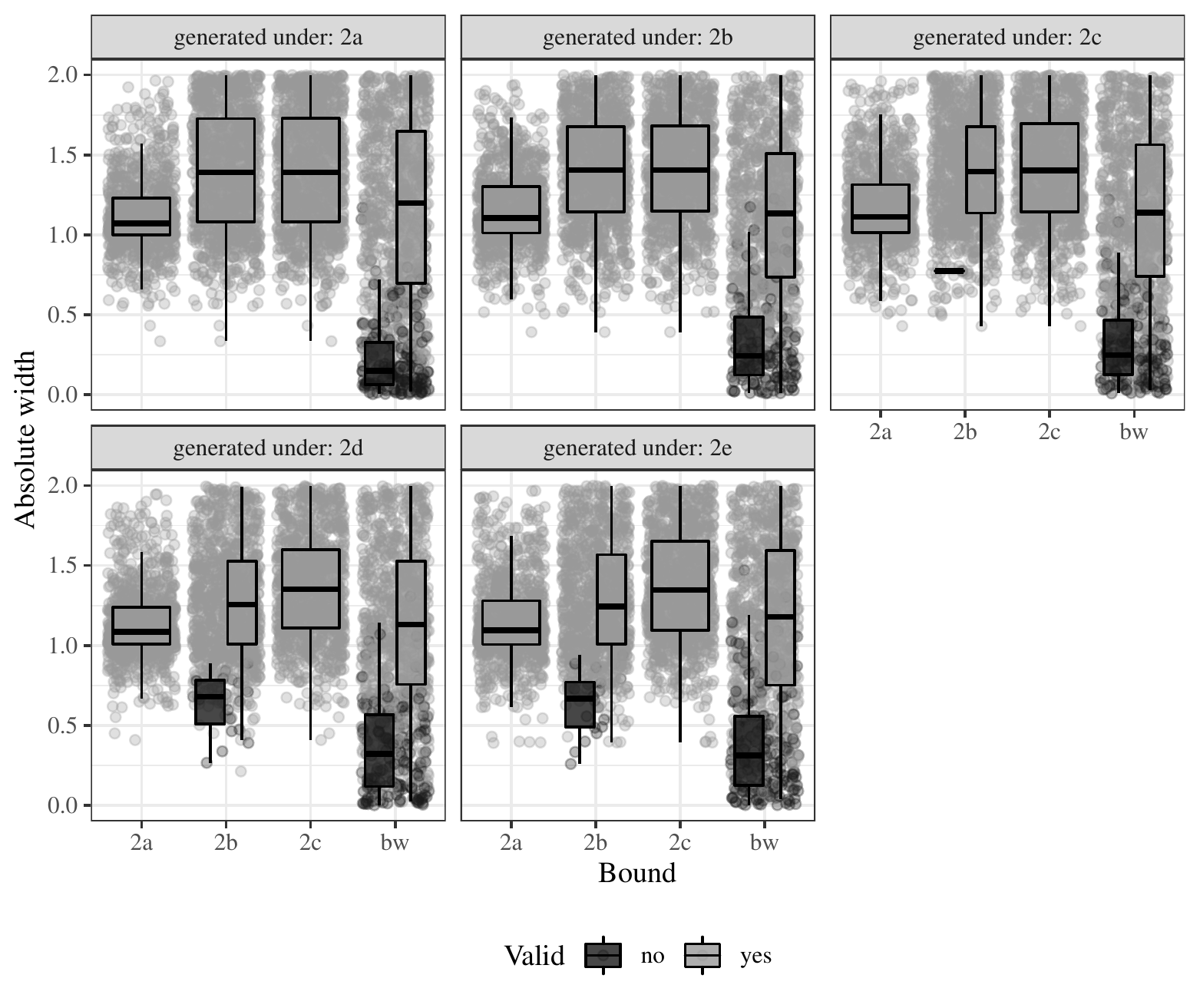}
\caption{\label{relwidth2} Comparison of the width of the true bounds for datasets generated under the DAGs in Figure 2 for distributions that are generated under Figure 2. The y-axis shows the absolute width of the bounds for each setting and the best/worst case bounds (denoted bw in the Figure). The light grey dots and boxes indicate cases where the bounds are valid (i.e., the true value is within the bounds), and the dark grey bounds that are invalid.}
\end{figure}

In Figure \ref{relwidth2}, we show the absolute width of the bounds for the settings of Figure 2. The bounds for the best/worst are frequently invalid for the intervention effect under all settings of Figure 2, as indicated by the darker shaded boxes and dots. The width of the bounds for the best/worst are clearly narrower, but since they target the assignment effect, $\tau$, they frequently do not cover the true intervention effect $\theta$. Under settings \ref{2c} - \ref{2e}, the bounds computed under \ref{2b} are occasionally invalid. The bounds of \ref{2a} seem quite robust, as we did not observe any distributions in which the bounds of \ref{2a} were invalid, this robustness was also seen in \citet{Gabriel2020}.

When generating distributions under $\alpha_3 > 0$, which is implied by the no defiers assumption, we find that the no defiers bounds for setting \ref{2b} are narrower than the \ref{2b} bounds allowing defiers 28\% of the time. The no defiers bounds for the other settings in Figure 2 are never narrower than the bounds allowing defiers for the same setting out of 10,000 generated  distributions. These results are illustrated in Figure S3 of supplementary materials.

\begin{figure}[ht]
\centering
\includegraphics[width = .95\textwidth]{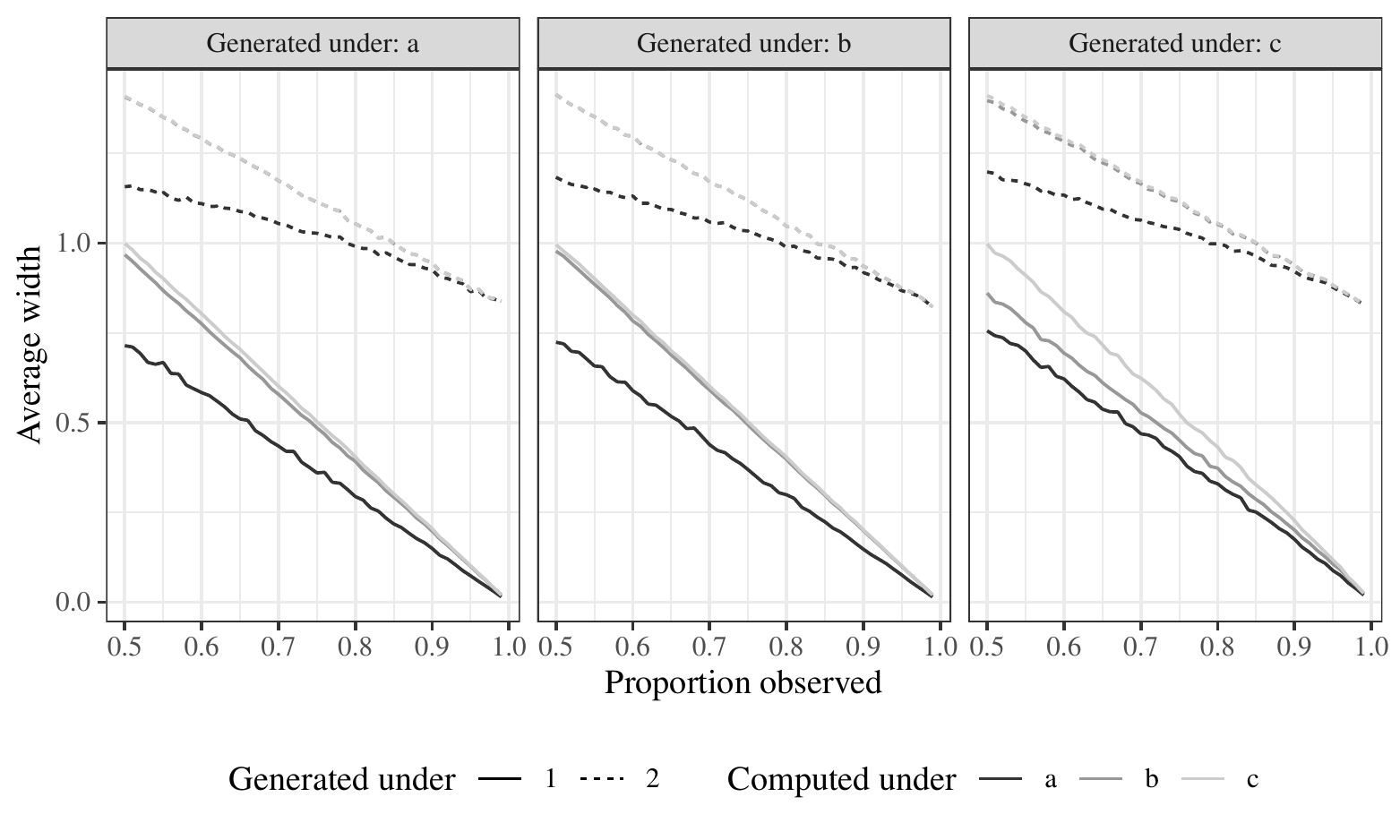}
\caption{\label{Fig:pobs} Illustration of the association between the proportion of observed outcomes and width of the bounds. The lines are the average width over the 1000 simulated distributions. Solid lines are under the models in Figure 1, dashed lines under Figure 2; dark grey is under model (a), medium grey under model (b), and light grey under model (c). }
\end{figure}

To investigate the impact of the amount of missingness on the informativeness of the study, we generate distributions with a fixed $\beta_3, \gamma_2$, and varied $\gamma_1$. Figure \ref{Fig:pobs} shows the average width of the bounds as functions of the proportion observed. Even with relatively small amounts of missing data $< 5\%$, we can see that the bounds quickly become very wide, particularly in the settings of Figure 2. The width of the bounds also appears to be approximately linearly increasing in the proportion missing.

To investigate the performance of the estimated bounds, we fix the values of the parameters and generate trials of size $n = 200$ or $2000$ from those distributions, calculate the empirical proportions needed to compute the bounds. We then use the nonparametric bootstrap of this procedure to compute quantile based 95\% confidence limits for the lower and upper bounds. Coverage of the 95\% bootstrap confidence intervals for the estimated bounds are shown in Table 1 for trial sizes of 200 and 2000 for a missingness probability of 25\%. We consider several values of $\theta$, over 1000 simulated replicates. We observe that a few of the confidence intervals have somewhat too small or too large coverage probability, but most have nearly 95\% coverage, as expected. Using the upper confidence limit of the upper bound and the lower confidence limit of the lower bound, we observe 100\% coverage of the true risk difference in these scenarios. 

\begin{table}
\begin{center}
\caption{Coverage of 95\% bootstrap CI for the true upper and lower bounds} 
\begin{tabular}{c|p{5.9mm}|c|c|c|c}
\hline
&&\multicolumn{4}{c}{Causal Diagrams}\\
true $\theta$ & trial size & 1a & 1b & 1c & 2c-2e\\
&&(Lower,Upper)&(Lower,Upper)&(Lower,Upper)&(Lower,Upper)\\
\hline
-0.2 & 200 & (0.92, 0.94) & (0.92, 0.94) & (0.92, 0.94) & (0.95, 0.94)\\
\hline
-0.2 & 2000 & (0.95, 0.95) & (0.95, 0.96) & (0.95, 0.96) & (0.94, 0.95)\\
\hline
-0.1 & 200 & (0.97, 0.95) & (0.96, 0.94) & (0.96, 0.94) & (0.92, 0.92)\\
\hline
-0.1 & 2000 & (0.97, 0.95) & (0.96, 0.94) & (0.96, 0.94) & (0.95, 0.93)\\
\hline
0.0 & 200 & (0.92, 0.96) & (0.92, 0.95) & (0.92, 0.95) & (0.94, 0.94)\\
\hline
0.0 & 2000 & (0.96, 0.96) & (0.95, 0.95) & (0.95, 0.95) & (0.95, 0.97)\\
\hline
\end{tabular}
\end{center}
\end{table}

\section{Real Data Application} \label{real}
\citet{du2015randomized} present the findings from a randomized controlled trial designed to estimate the causal effect of peanut consumption on the development of allergy to peanuts in infants. 640 participants between 4 months and 11 months of age were randomized to either consume peanuts or avoid peanuts until the age of 60 months. Compliance with the assigned intervention was assessed weekly by using a food frequency questionnaire, and by manual inspection of the infants' cribs for peanut crumbs in a subset of participants. At 60 months, the primary outcome of peanut allergy was assessed using an oral food challenge. Outcome data were missing in some participants either due to loss to follow up, or due to failure of the oral food challenge procedures. The publicly available trial data were downloaded from the The Immune Tolerance Network TrialShare website on 2020-06-15 (\url{https://www.itntrialshare.org/}, study identifier: ITN032AD). %R code to reproduce the bounds we report is available as supplementary material. 

This study clearly falls into one of the settings of Figure 2, as both compliance and missing outcome data were issues in the study. The primary results in the manuscript were reported as the proportion with food allergy at 60 months in the assigned intervention groups. The per-protocol analysis and the worst case imputation analysis were reported as sensitivity analyses. Here we compute and report our bounds.  

Our estimated bounds for $\theta$ and $\tau$ are shown in Table 2 along with bootstrap 95\% confidence intervals. We see that noncompliance and missing data lead to a great deal more uncertainty in the causal effect estimate relative to sampling variability. Nevertheless, the bounds still exclude the risk difference of 0, suggesting that there is compelling evidence that consuming peanuts reduces the risk of peanut allergy at 60 months. Compared to the point estimate of $-0.14$ reported by \citet{du2015randomized}, the range of possible causal effects goes from $-0.01$ to $-0.29$ without any additional assumptions. The original publication reports the per-protocol estimate of the intervention effect as $-0.17$, and the worst case imputation estimate as $-0.12$. Based on inspection of the publicly available data, however, their worst case imputation estimate is more accurately described as a ``pessimistic imputation'', rather than worst case, since not all subjects missing outcomes in the intervention arm were imputed with having an allergic event and not all subjects missing outcomes in the avoidance arm were imputed as not having an event. Thus, our best/worst case bounds cover, but are not exactly the same as their published ``worst case'' imputation results.

\begin{table}
\begin{center}
\caption{Bounds with 95\% confidence intervals in the peanut allergy trial}
\begin{tabular}{l|ll}
%\hline
 & \multicolumn{2}{c}{Intervention effect $\theta$}  \\
 & Lower bound & Upper bound \\
\hline
2a & $\cithree{-0.29}{-0.25}{-0.19}$ & $\cithree{-0.11}{-0.06}{-0.01}$ \\
\hline
2b & $\cithree{-0.29}{-0.25}{-0.20}$ & $\cithree{-0.12}{-0.06}{-0.02}$ \\
2b no defiers & $\cithree{-0.29}{-0.25}{-0.20}$ & $\cithree{-0.15}{-0.09}{-0.04}$\\
\hline
2c-2e & $\cithree{-0.29}{-0.25}{-0.20}$ & $\cithree{-0.11}{-0.06}{-0.01}$ \\
\hline
& \multicolumn{2}{c}{Assignment effect $\tau$} \\
& Lower bound & Upper bound \\
\hline
best/worst & $\cithree{-0.22}{-0.17}{-0.12}$ & $\cithree{-0.14}{-0.10}{-0.06}$\\
1a & $\cithree{-0.21}{-0.16}{-0.11}$ & $\cithree{-0.18}{-0.14}{-0.10}$ \\
1b & $\cithree{-0.22}{-0.17}{-0.12}$ & $\cithree{-0.16}{-0.10}{-0.06}$\\
1c & $\cithree{-0.22}{-0.17}{-0.12}$ & $\cithree{-0.14}{-0.10}{-0.06}$ 
\end{tabular}
\end{center}
\end{table}

\section{Discussion} \label{dis}
To ensure validity of causal effect estimates in a randomized experiment, every effort should be made to avoid missing data due to drop-out  \citep{fleming2011addressing}. When missing data are unavoidable, our bounds can be used to quantify the uncertainty in the causal effect of an intervention while making minimal assumptions about the nature of the missingness mechanism. Our bounds can often be narrower than the best/worst case bounds in settings with perfect compliance. It is also of note that although the technique of best/worst sensitivity analysis is commonly applied and reported in clinical trials, to our knowledge the nonparametric bounds implied by the procedure based on the true probabilities have not been previously presented in this manner in the literature. 

When noncompliance is also an issue our proposed bounds provide direct information on the causal effect of the intervention, in contrast to the best/worst case imputation approach which assesses the effect of assignment to intervention. Additionally, when no defiers is a plausible assumption, our bounds can be tightened in particular settings. Our motivating data example demonstrates how our bounds can be applied to answer important scientific questions regarding the size of causal effects in trials that are subject to noncompliance and nonignorable missing data. 

We have assumed throughout that the practitioner, having randomized the experiment and followed its progression, has adequate knowledge to determine the underlying causal diagram. We acknowledge that this may not be the case. It may, however, be possible to use the observed data in some settings to infer the causal relationships via causal discovery algorithms \citep{spirtes1991algorithm}, or by observing that the computed bounds are not compatible with the assumed settings, i.e. the computed upper bound is less than the computed lower bound. This is of course a limitation of this work as in settings where the assumed causal diagram does not hold, the bounds are in no way guaranteed to cover the true causal effect. However, unlike observational settings, there are many characteristics of the experiment that can help narrow the set of plausible causal diagrams without the need to test. For example, in triple blind clinical trials it is implausible that randomization would have a causal effect on missingness, or that in a point-of-care single time point intervention there would be noncompliance. These characteristics should clearly be considered when selecting the assumed setting under which to calculate the bounds. 

Although we have considered the addition of the no defiers assumption in settings with noncompliance, there are many additional monotonicity assumptions that could be made in the various settings. For example, it may be plausible that missingness is monotone in the intervention or outcome in some settings, which may lead to tighter bounds. Additionally, the stronger assumption that no control subject, $R=0$, can take the intervention $X=1$, may lead to tighter or simply different bounds than have been derived under the weaker no defiers assumption. Investigation of such additional monotonicity settings is a current area of research for the authors. 

\bibliographystyle{abbrvnat}
\bibliography{refer.bib}

\end{document}

\def\spacingset#1{\renewcommand{\baselinestretch}%
{#1}\small\normalsize} \spacingset{0}

%%%%%%%%%%%%%%%%%%%%%%%%%%%%%%%%%%%%%%%%%%%%%%%%%%%%%%%%%%%%%%%%%%%%%%%%%%%%%%

\if1\blind
{
  \title{\bf Supplementary materials for ``Nonparametric bounds for causal effects in imperfect randomized experiments"}
  \author{Erin E. Gabriel$^{1}$\thanks{EEG is partially supported by Swedish Research Council grant 2017-01898, AS by Swedish Research Council grant 2016-01267 and MCS by Swedish Research Council grant 2019-00227} \hspace{.2cm}\\
  Arvid Sj\"olander$^1$\\
  Michael C. Sachs$^1$\\
\footnotesize 1:Department of Medical Epidemiology and Biostatistics, Karolinska Institutet, Stockholm, Sweden \\
$\star$The authors contributed equally to this work.\\
corresponding author: erin.gabriel@ki.se}
  \maketitle
} \fi 

\if0\blind
{
  \bigskip
  \bigskip
  \bigskip
  \begin{center}
    {\LARGE\bf Supplementary materials for ``Nonparametric bounds for causal effects in imperfect randomized experiments"}
\end{center}
  \medskip
} \fi

\bigskip

\spacingset{2}
\vspace{-1cm}
\section{Previous Bounds}

For convenience of notation, we define the following probability abbreviations. Let 
\begin{eqnarray*}
p_{y.x1} &=& p\{Y=y|X=x,O=1\}, \\
p_{y1.x} &=& p\{Y=y, O=1|X=x\}, \\
p_{xy1.r} &=& p\{X=x,Y=y,O=1|R=r\},\\
p_{xy.1r} &=& p\{X=x,Y=y|O=1,R=r\},\\
p_{o.x} &=& p\{O=o|X=x\},\\
q_{xy.o} &=& p\{X=x, Y=y|O=o\},\\
q_{r} &=& p\{O=1|R=r\}, \\
q &=& p\{O=1\}.
\end{eqnarray*}

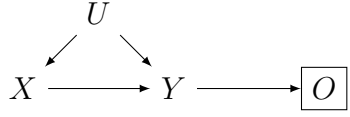
\begin{figure}[ht]
\begin{center}
\begin{tikzpicture}
\node (uh) at (1,1) {$U$};
\node (v) at (0,0) {$X$};
\node (i) at (2,0) {$Y$};
\node[draw, rectangle] (s) at (4,0) {$O$};
\draw[-latex] (v) -- (i);
\draw[-latex] (i) -- (s);

\draw[-latex] (uh) -- (v);
\draw[-latex] (uh) -- (i);

%\draw[-latex] (uh) -- +(2.79,-0.9) (i);

\end{tikzpicture}
\end{center}
\caption{\label{supfig} Setting with unconfounded missingness in the outcome.}
\end{figure}

In \citet{Gabriel2020}, bounds were derived for $\theta$ under the causal diagram in Figure \ref{supfig}. 

Define $A(y,o)=q_{1y.o}/q_{0y.o}$. The bounds derived by \citet{Gabriel2020} are
\begin{eqnarray}
\label{Beq:5}
\theta\geq -(q_{10.1}+q_{01.1})q-\max\left\{\frac{1}{1+A(1,1)},\frac{A(0,1)}{1+A(0,1)}\right\}(1-q)
\end{eqnarray}
and
\begin{eqnarray}
\label{Beq:6}
\theta\leq 1-(q_{10.1}+q_{01.1})q-\min\left\{\frac{1}{1+A(1,1)},\frac{A(0,1)}{1+A(0,1)}\right\}(1-q).
\end{eqnarray}
We will use these bounds in the refinement of our proposed Figure 1a bounds.

Bounds derived in \citet{Gabriel2020} for the setting of Figure 2a from the main text without the no defiers assumption, are provided below in terms of the abbreviations that match the main text, where  $B(y,r,o)=p_{1y.ro}/p_{0y.ro}$.
\begin{eqnarray}
\label{Beq:S1}
\theta\geq \max\left\{\begin{array}{l}
p_{11.11}q_1+p_{00.01}q_0-1\\
p_{11.10}q_0+p_{00.11}q_1-1\\
(p_{11.10}-p_{10.10}-p_{01.10})q_0-(p_{11.11}+p_{01.11})q_1\\ \phantom{11111} -\frac{B(0,0,1)}{1+B(0,0,1)}(1-q_0)-(1-q_1)\\
(p_{11.11}-p_{10.11}-p_{01.11})q_1-(p_{11.10}+p_{01.10})q_0\\
\phantom{11111}-\frac{B(0,1,1)}{1+B(0,1,1)}(1-q_1)-(1-q_0)\\
-(p_{10.11}+p_{01.11})q_1-\max\left\{\frac{1}{1+B(1,1,1)},\frac{B(0,1,1)}{1+B(0,1,1)}\right\}(1-q_1)\\
 -(p_{10.10}+p_{01.10})q_0-\max\left\{\frac{1}{1+B(1,0,1)},\frac{B(0,0,1)}{1+B(0,0,1)}\right\}(1-q_0)\\
 (p_{00.11}-p_{10.11}-p_{01.11})q_1-(p_{10.10}+p_{00.10})q_0\\
\phantom{11111}-\max\left\{\frac{1}{1+B(1,1,1)},-\frac{1-B(0,1,1)}{1+B(0,1,1)}\right\}(1-q_1)-(1-q_0)\\
(p_{00.10}-p_{10.10}-p_{01.10})q_0-(p_{10.11}+p_{00.11})q_1\\
\phantom{11111}-\max\left\{\frac{1}{1+B(1,0,1)},-\frac{1-B(0,0,1)}{1+B(0,0,1)}\right\}(1-q_0)-(1-q_1)
\end{array}\right\}, 
\end{eqnarray}
and 
\begin{eqnarray}
\label{Beq:S2}
\theta\leq \min\left\{\begin{array}{l}1-p_{10.11}q_1-p_{01.10}q_0\\
1-p_{10.10}q_0-p_{01.11}q_1\\
(p_{10.11}+p_{00.11})q_1+(p_{11.10}+p_{00.10}-p_{10.10})q_0\\
\phantom{11111}+(1-q_1)+\max\left\{\frac{1-B(0,0,1)}{1+B(0,0,1)},\frac{B(1,0,1)}{1+B(1,0,1)}\right\}q_0\\
(p_{11.11}+p_{00.11}-p_{10.11})q_1+(p_{10.10}+p_{00.10})q_0\\
\phantom{11111}+\max\left\{\frac{1-B(0,1,1)}{1+B(0,1,1)},\frac{B(1,1,1)}{B(1,1,1)}\right\}(1-q_1)+(1-q_0)\\
(p_{11.11}+p_{00.11})q_1+\max\left\{\frac{1}{1+B(0,1,1)},\frac{B(1,1,1)}{1+B(1,1,1)}\right\}(1-q_1)\\
(p_{11.10}+p_{00.10})q_0+\max\left\{\frac{1}{1+B(0,0,1)},\frac{B(1,0,1)}{1+B(1,0,1)}\right\}(1-q_0)\\
(p_{11.11}+p_{00.11}-p_{01.11})q_1+(p_{11.10}+p_{01.10})q_0\\
\phantom{11111}+\max\left\{\frac{1}{B(0,1,1)},-\frac{1-B(1,1,1)}{1+B(1,1,1)}\right\}(1-q_1)+(1-q_0)\\
(p_{11.10}+p_{00.10}-p_{01.10})q_0+(p_{11.11}+p_{01.11})q_1\\
\phantom{11111}+\max\left\{\frac{1}{B(0,0,1)},-\frac{1-B(1,0,1)}{1+B(1,0,1)}\right\}(1-q_0)+(1-q_1)
\end{array}\right\}.
\end{eqnarray}

Bounds derived in \citet{Gabriel2020} for the setting of Figure 2b from the main text without the no defiers assumption, are provided below in terms of the abbreviations that match our main text. 
\begin{eqnarray}
\label{Beq:11}
\theta\geq \max\left\{\begin{array}{l}
p_{001.1} - p_{011.0} - p_{111.0} + 2p_{111.1} - 1\\
p_{001.1} + p_{111.0} - 1\\
p_{001.0} + p_{111.1} - 1\\
p_{001.1} + p_{111.1} - 1\\
p_{001.0} - p_{011.1} + 2p_{111.0} - p_{111.1} - 1\\
- p_{001.0} + 2p_{001.1} - p_{101.0} + p_{111.1} - 1\\
p_{001.0} + p_{111.0} - 1\\
2p_{001.0} - p_{001.1} - p_{101.1} + p_{111.0} - 1\\
2p_{001.0} - p_{001.1} - p_{101.1} - p_{011.1} + 2p_{111.0} - p_{111.1} - 1\\
- p_{001.0} + 2p_{001.1} - p_{101.0} - p_{011.0} - p_{111.0} + 2p_{111.1} - 1
\end{array}\right\}.
\end{eqnarray}

\begin{eqnarray}
\label{Beq:12}
\theta\leq \min\left\{\begin{array}{l}
- p_{101.0} - 2p_{011.0} + p_{011.1} + p_{111.1} + 1\\
- p_{101.1} + p_{011.0} - 2p_{011.1} + p_{111.0} + 1\\
p_{001.1} - 2p_{101.0} + p_{101.1} - p_{011.0} + 1\\
- p_{101.0} - p_{011.1} + 1\\
- p_{101.1} - p_{011.0} + 1\\
- p_{101.0} - p_{011.0} + 1\\
p_{001.1} - 2p_{101.0} + p_{101.1} - 2p_{011.0} + p_{011.1} + p_{111.1} + 1\\
- p_{101.1} - p_{011.1} + 1\\
p_{001.0} + p_{101.0} - 2p_{101.1} - p_{011.1} + 1\\
p_{001.0} + p_{101.0} - 2p_{101.1} + p_{011.0} - 2p_{011.1} + p_{111.0} + 1\\
\end{array}\right\}.
\end{eqnarray}

\section{Linear programming method}
In order to use this algorithm one must show that the problem in question can be stated as a linear programming problem. First, we observe that the causal diagrams in Figure 1 are equivalent to the diagrams in which all latent parents of the variables $X, Y, O$ are replaced with the response function variables $W_X, W_Y, W_O$ which are categorical. Since $Y$ and $O$ have common cause $U$ in the original diagrams, their response function variables $W_Y, W_O$ are dependent but independent from $W_X$. Further, since $X, Y, O$ are all binary, the response function variables are each categorical with $2^{k_i}$ levels, where $k_i$ is the number of parents of variable $i (\in \{X, Y, O\})$ in the original diagram. The response function variables specify the manner in which the variables are determined from their parents in the diagram. Specifically, each variable is functionally determined from the value of its response function variable and the value of its parents; i.e, there exist functions $f_O$, $f_Y$, and $f_X$ such that

\begin{eqnarray*}
x &=& f_X(\mathbf{w}) = f_X(w_x)  \\
y &=& f_Y(\mathbf{w}) = f_Y(w_y, x) \\
o &=& f_O(\mathbf{w}) = f_O(w_o, y, x)
\end{eqnarray*}
are defined for each possible level of $\mathbf{w} = (w_x, w_y, w_o)$. This representation permits us to relate potential outcome probabilities to probabilities of the response function variables. For example, $p\{Y(X = x) = 1\} = p\{f_Y(w_y, x) = 1\} = \sum_{w_y \in \Gamma(x)}p\{W_y = w_y\}$ where $\Gamma(x) = \{w_y: f_Y(w_y, x) = 1\}$.

To determine the constraints implied by the causal diagram, we can relate all observed probabilities to the probabilities of the response function variables by recursively evaluating the response functions, e.g., for a fixed $\mathbf{w}$, the observed variables are determined by 

\begin{eqnarray*}
x &=& f_X(w_x) \\
y &=& f_Y(w_y, f_X(w_x)) \\
o &=& f_O(w_o, f_Y(w_y, f_X(w_x)), f_X(w_x)). 
\end{eqnarray*}

Therefore, for $x, y \in \{0, 1\}$
$$p\{X = x, Y = y, O = 1\} = \sum_{\mathcal{W}} p\{W_X = w_x, W_Y = w_y, W_O = w_o\}, $$
where 
$$
\mathcal{W} = \{\mathbf{w}: f_X(\mathbf{w}) = x, f_Y(\mathbf{w}) = y, f_O(\mathbf{w}) = 1\}.
$$
Since $W_X$ is independent of $(W_Y, W_O)$, we can factorize the terms in the sum

\begin{eqnarray*}
p\{X = x, Y = y, O = 1\} & = & p\{W_X = w_x\} \sum_{\mathcal{W}^-}  p\{W_Y = w_y, W_O = w_o\} \implies \\
p\{Y = y, O = 1 | X = x\} & = & \sum_{\mathcal{W}^-}  p\{W_Y = w_y, W_O = w_o\}, 
\end{eqnarray*}
where $\mathcal{W}^- = \{(w_x, w_y): f_X(\mathbf{w}) = x, f_Y(\mathbf{w}) = y\}$. The second equation follows because $p\{W_x = w_x\} = p\{X = x\}$ for $x \in \{0, 1\}$ in the diagrams of Figure 1. The probabilities $p\{O = 0, X = x\}$, $x \in \{0, 1\}$, which may be known, also define constraints, but it can easily be shown that these constraints are redundant and can be eliminated. Therefore, the constraints on the conditional probabilities are linear in the response function variables, necessary, and sufficient for the distribution to be in the causal model for Figures 1b and 1c of the main text. This, combined with the fact that the risk difference is linear in the response function variable probabilities means that the bounds on the risk difference can be obtained as the solution to the linear programming problem. This line of argument applies for Figures 1b and 1c of the main text, with the only difference being whether $f_O$ depends on $x$.

For Figures 2b-2e of the main text, the derivations become slightly more complex due to the additional variable. Nevertheless, the same logic allows us to factorize the response function variable for $R$ out of the constraints, and therefore derive necessary and sufficient constraints on the probability distribution that are linear, and in terms of probabilities of the form $p_{xy1.r}$. 

The linear objective (the risk difference) and linear constraints defines a constrained linear optimization problem. Solutions to this problem can be found symbolically by applying Balke's implementation of a vertex enumeration algorithm \citep{balke1994counterfactual, mattheiss1973algorithm}. In brief, this algebraically reduces the variables in the optimization problem, then adds slack variables so that all constraints are converted into inequality constraints. The dual of this problem is to maximize (minimize) a linear function of the observed probabilities, subject to a set of constraints. Thus the extremum of the causal query as stated in terms of the potential outcome probabilities is equal to the extremum in the observed probability space defined by the dual constraints. Then, by noting that those constraints describe a convex polytope in the observed probability space, the global extrema can be found by enumerating all of the vertices of the polytope. This gives the bounds on the causal effect of interest as the minimum (maximum) of a list of terms involving only observable probabilities, each of which corresponds to a vertex of this polytope. This demonstrates that for these problems in Figures 1b, 1c, and 2b-2e, tight and valid bounds on $\theta$ can be derived symbolically according to this algorithm.

\section{Derivation of Theorem 1}
In the setting of Figure 1a of the main text, the data distribution is defined by the 8 probabilities $q_{xy.o}$ for $(x,y,o) \in{0,1}$. The unknown probabilities $q_{xy.0}, (x,y)\in{0,1}$, are constrained by
\begin{eqnarray}
\label{eq:1}
\sum_{xy}q_{xy.0}=1.
\end{eqnarray} 
If we assume that $0<p\{O=1|Y=y\}<1$, and that all observed probabilities are $0<p_{y.x1}<1$, then  
\begin{eqnarray}
\label{eq:Ay}
A(y,0)=A(y,1)\textrm{ for all } y.
\end{eqnarray}
This can be seen because $O\bot X|Y$ and thus 
\begin{eqnarray*}
\frac{p\{X=1|Y=y\}}{p\{X=0|Y=y\}}
&=&\frac{p\{X=1|Y=y,O=o\}}{p\{X=0|Y=y,O=o\}}\\
&&\phantom{1}=\frac{q_{1y.o}/p\{Y=y|O=o\}}{q_{0y.o}/p\{Y=y|O=o\}}\\
&&\phantom{1}=A(y,o)
\end{eqnarray*}
does not depend on the value of $o$. 

Starting from the identified point 
\begin{eqnarray}
\label{Beq:1}
\theta = 1-(p\{Y=0|X=1\}+p\{Y=1|X=0\})
\end{eqnarray}  
we can expand and derive bounds in our setting with partially observed data. To derive bounds for $\theta$ in the setting of non-ignorable missingness, Figure 1a, that take the constraint in \eqref{eq:Ay} into account we proceed as follows. We partition \ref{Beq:1} into observed and unobserved pieces.
\begin{eqnarray}
\label{eq:partitioning}
1-(p\{Y=0|X=1\}+p\{Y=1|X=0\})=\nonumber \\ 1-\left(\frac{q_{10.1}*p_{1.1}}{p\{X=1|O=1\}}+\frac{q_{01.1}*p_{0.1}}{p\{X=0|O=1\}}\right) \nonumber\\
-\left(\frac{q_{10.0}*p_{1.0}}{p\{X=1|O=0\}}+\frac{q_{01.0}*p_{0.0}}{p\{X=0|O=0\}}\right).
\end{eqnarray}
Under \eqref{eq:1} and \eqref{eq:Ay} we can express $(q_{01.0},q_{10.0},q_{11.0})$ as functions of $q_{00.0}$. From \eqref{eq:Ay} we have  
\begin{eqnarray}
\label{eq:4}
q_{10.0}=q_{00.0}A(0,1).
\end{eqnarray}
We further have 
\begin{eqnarray}
\label{eq:5}
q_{01.0}&=&1-q_{00.0}-q_{10.0}-q_{11.0}\nonumber\\
&=&1-q_{00.0}\{1+A(0,1)\}-q_{11.0},
\end{eqnarray}
where the first equality follows from (\ref{eq:1}) and the second from (\ref{eq:4}). Plugging (\ref{eq:5}) into (\ref{eq:Ay}) gives
\begin{eqnarray*}
\frac{q_{11.0}}{1-q_{00.0}\{1+A(0,1)\}-q_{11.0}}=A(1,1),
\end{eqnarray*}
which gives 
\begin{eqnarray}
\label{eq:6}
q_{11.0}=\frac{A(1,1)[1-q_{00.0}\{1+A(0,1)\}]}{1+A(1,1)}.
\end{eqnarray}
Finally, plugging (\ref{eq:6}) into (\ref{eq:5}) gives
\begin{eqnarray}
\label{eq:7}
q_{01.0}=\frac{1-q_{00.0}\{1+A(0,1)\}}{1+A(1,1)}.
\end{eqnarray}
We proceed by finding the feasible space for $q_{00.0}$, i.e. the space for which $(q_{00.0},q_{01.0},q_{10.0},q_{11.0})$ are all $\in [0,1]$. $(q_{01.0},q_{10.0},q_{11.0})$ are monotone functions of $q_{00.0}$. Thus, the feasible space for $q_{00.0}$ is a continuous range with $q_{00.0} > 0$. To find the upper bound for $q_{00.0}$, note that $q_{10.0}$ is increasing in $q_{00.0}$, with $q_{10.0}$ approaching 1 when $q_{00.0}$ approaches $1/A(0,1)$; $q_{01.0}$ is decreasing in $q_{00.0}$, with $q_{01.0}$ approaching zero when $q_{00.0}$ approaches $1/\{1+A(0,1)\}$; $q_{11.0}$ is decreasing in $q_{00. 0}$, with $q_{11.0}$ approaching zero when $q_{00.0}$ approaches $1/\{1+A(0,1)\}$. Thus,  $q_{00.0}\leq \min[1/A(0,1),1/\{1+A(0,1)\}]=1/\{1+A(0,1)\}$. We thus have that 
\begin{eqnarray}
\label{eq:8}
0\leq q_{00.0} \leq \frac{1}{1+A(0,1)}.
\end{eqnarray}
From (\ref{eq:4}) and (\ref{eq:7}) we have that 
\begin{eqnarray}
\label{eq:9}
q_{01.0}+q_{10.0}=\frac{1+q_{00.0}[A(0,1)\{1+A(1,1)\}-\{1+A(0,1)\}]}{1+A(1,1)}.
\end{eqnarray}
If $A(0,1)/\{1+A(0,1)\}\geq 1/\{1+A(1,1)\}$, then the right hand side of (\ref{eq:9}) is increasing in $q_{00.0}$. Using (\ref{eq:8}) we thus have that 
\begin{eqnarray}
\label{eq:14}
\min(q_{01.0}+q_{10.0})&=&\frac{1}{1+A(1,1)},\nonumber\\
\max(q_{01.0}+q_{10.0})&=&\frac{A(0,1)}{1+A(0,1)}.
\end{eqnarray} 
If $A(0,1)/\{1+A(0,1)\}\leq 1/\{1+A(1,1)\}$, then the right hand side of (\ref{eq:9}) is decreasing in $q_{00.0}$. Using (\ref{eq:8}) we thus have that 
\begin{eqnarray}
\label{eq:15}
\min(q_{01.0}+q_{10.0})&=&\frac{A(0,1)}{1+A(0,1)},\nonumber\\
\max(q_{01.0}+q_{10.0})&=&\frac{1}{1+A(1,1)}.
\end{eqnarray}
Combining (\ref{eq:14}) and (\ref{eq:15}) gives
\begin{eqnarray}
\label{eq:16}
\min(q_{01.0}+q_{10.0})&=&\min\left\{\frac{1}{1+A(1,1)},\frac{A(0,1)}{1+A(0,1)}\right\},\nonumber\\
\max(q_{01.0}+q_{10.0})&=&\max\left\{\frac{1}{1+A(1,1)},\frac{A(0,1)}{1+A(0,1)}\right\}.
\end{eqnarray}

\begin{small}
\begin{eqnarray}
\label{eq:17}
\min\left(\frac{q_{10.0}*p_{1.0}}{p\{X=1|O=0\}}+\frac{q_{01.0}*p_{0.0}}{p\{X=0|O=0\}}\right) \geq  \min\left\{\frac{1}{1+A(1,1)},\frac{A(0,1)}{1+A(0,1)}\right\}* \nonumber \\ \min\left\{\frac{p_{1.0}}{p\{X=1|O=0\}}, \frac{p_{0.0}}{p\{X=0|O=0\}}\right\}, \nonumber \\
\max\left(\frac{q_{10.0}*p_{1.0}}{p\{X=1|O=0\}}+\frac{q_{01.0}*p_{0.0}}{p\{X=0|O=0\}}\right) \leq \max\left\{\frac{1}{1+A(1,1)},\frac{A(0,1)}{1+A(0,1)}\right\}* \nonumber \\ \max\left\{\frac{p_{1.0}}{p\{X=1|O=0\}}, \frac{p_{0.0}}{p\{X=0|O=0\}}\right\}  . 
\end{eqnarray}
\end{small}

Combining (\ref{eq:partitioning}) and (\ref{eq:17}) gives
\begin{small}
\begin{eqnarray}
\label{Beq:18}
&&\theta\geq 1-\left(p_{0.11}*p_{1.1}+p_{1.01}*p_{1.0}\right)-\max\left\{\frac{1}{1+A(1,1)},\frac{A(0,1)}{1+A(0,1)}\right\} \nonumber \\ 
&&*\max\left\{\frac{p_{0.1}}{p\{X=1|O=0\}},\frac{p_{0.0}}{p\{X=0|O=0\}}\right\} 
\end{eqnarray}
and 
\begin{eqnarray}
\label{Beq:19}
&&\theta\leq 1-\left(p_{0.11}*p_{1.1}+p_{1.01}*p_{1.0}\right)-\min\left\{\frac{1}{1+A(1,1)},\frac{A(0,1)}{1+A(0,1)}\right\} \nonumber \\ 
&&*\min\left\{\frac{p_{0.1}}{p\{X=1|O=0\}},\frac{p_{0.0}}{p\{X=0|O=0\}}\right\} 
\end{eqnarray}
\end{small}

Therefore the bounds given in the main text for $\theta$ are valid in the setting of Figure 1a, provided that missingness was not deterministic, i.e. $0<p\{O=1|Y=y\}<1$  and that all observed probabilities are $0<p_{y.x1}<1$. 

If $A(0,1)$ or $A(1,1)$ are undefined but $A(0,1)^{-1}$ or $A(1,1)^{-1}$ are not undefined, then switching the coding of the randomization $x*=1-x$, then the bounds for $l^*\leq\theta^*\leq u^*$, can be used to bound $\theta=-\theta^*$ by $-u^* \leq\theta\leq -l^*$. If both $A(0,1)$ and $A(0,1)^{-1}$ or $A(0,1)$ and $A(0,1)^{-1}$ are undefined, then the best that can be done is to bound the max of the unobserved and undefined quantities by 1, and the min by 0. Since we will be able to observe if $0<p\{O=1|Y=y\}$, but not if $p\{O=1|Y=y\}<1$, as we will have no way of knowing if we have observed all the subjects with the outcome $y$ or if there are subjects with the event or outcome that are missing. 

However, provided $0<p\{O=1|Y=y\}$ then if $A(1,1)$ or $A(1,1)^{-1}$ are bounded then the bounds will be valid, but not as tight as they could be if one knew that $p\{O=1|Y=y\}=1$ and could set one or more of unobserved probabilities to zero. When  $p\{O=1|Y=y\}=0$ for any value of $y$, this will be clearly observable, and $A(y,1)$ or $A(y,1)^{-1}$ will both be undefined. Although it is possible to bound the causal effect by replacing the max of unknown probabilities with 1 and the min of unknown probability with 0, this will never be truly informative and instead we suggest evaluating what went wrong in the trial and either abandoning the intervention or re-running the study altogether. 

For example, if one observes no cases implying that $p\{O=1|Y=1\}=0$, this may be because the event is extremely rare, the sample size is too small or follow-up time too short. Alternatively if one observes that $p\{O=1|Y=0\}=0$, either the intervention does not work and the event is extremely common or the follow-up time is too long and the outcome is inevitable, e.g. death or taxes. Additionally, if $p\{X=x|O=0\}=0$ for any $x$, such that the ratios are undefined, this may be due to the intervention having a direct effect on sampling, suggesting that the causal diagram in Figure 1a of the main text may not be the correct setting.

\section{Comparison and refinement of bounds} 
\begin{figure}[ht]
\centering
\includegraphics[width = .95\textwidth]{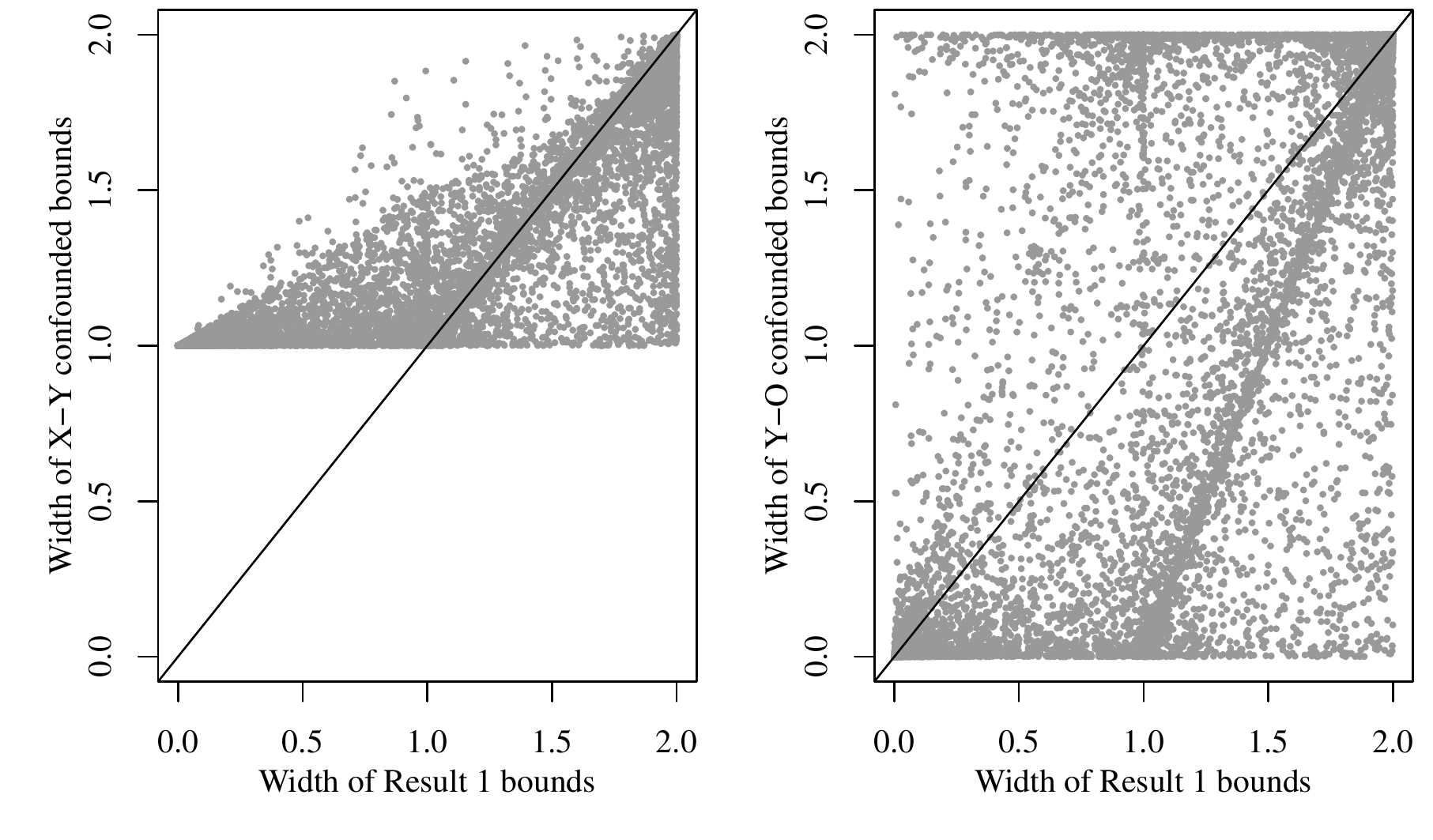}
\caption{Comparison of the width of the bounds for Theorem 1 to the bounds under one of the types of confounding. All of the bounds are valid, yet sometimes the bounds under confounding are narrower than those derived in Theorem 1. \label{ref1c}}
\end{figure}

\clearpage

\section{Additional simulation Results}

\begin{figure}[ht]
\centering
\includegraphics[width = .99\textwidth]{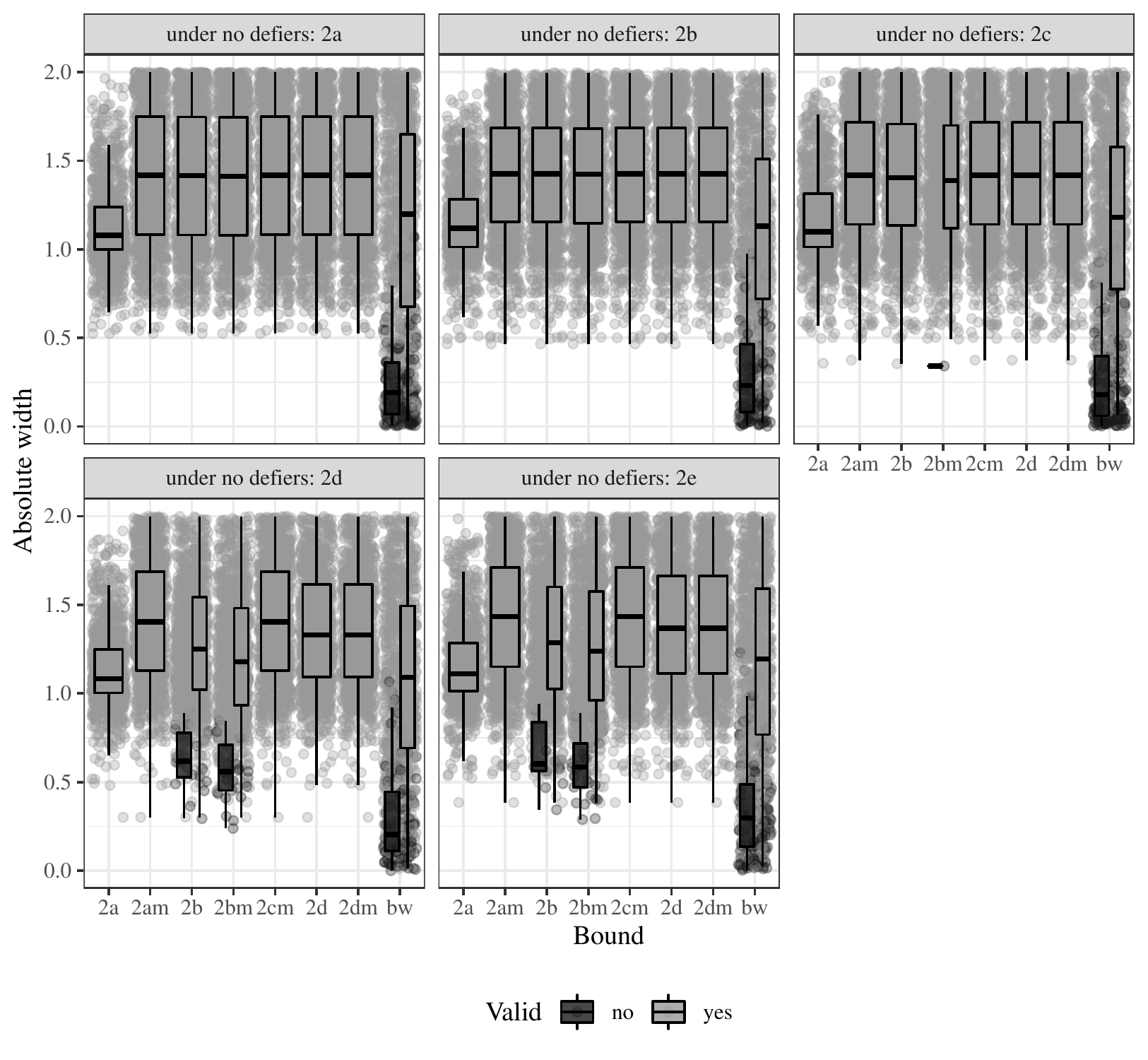}
\caption{\label{monotonewidth} Comparison of the width of the true bounds for the DAGs in Figure 2 for distributions that are generated under $\alpha_3 > 0$ which is implied by the no defiers assumption. The bounds ending in ``m'' are the ones derived under the no defiers assumption. The light grey dots and boxes indicate cases where the true value is within the bounds (valid), and the dary grey where the true value is not within the bounds (invalid).}
\end{figure}

\begin{figure}[ht]
\begin{subfigure}{.95\textwidth}
\centering
\includegraphics[width = .99\textwidth]{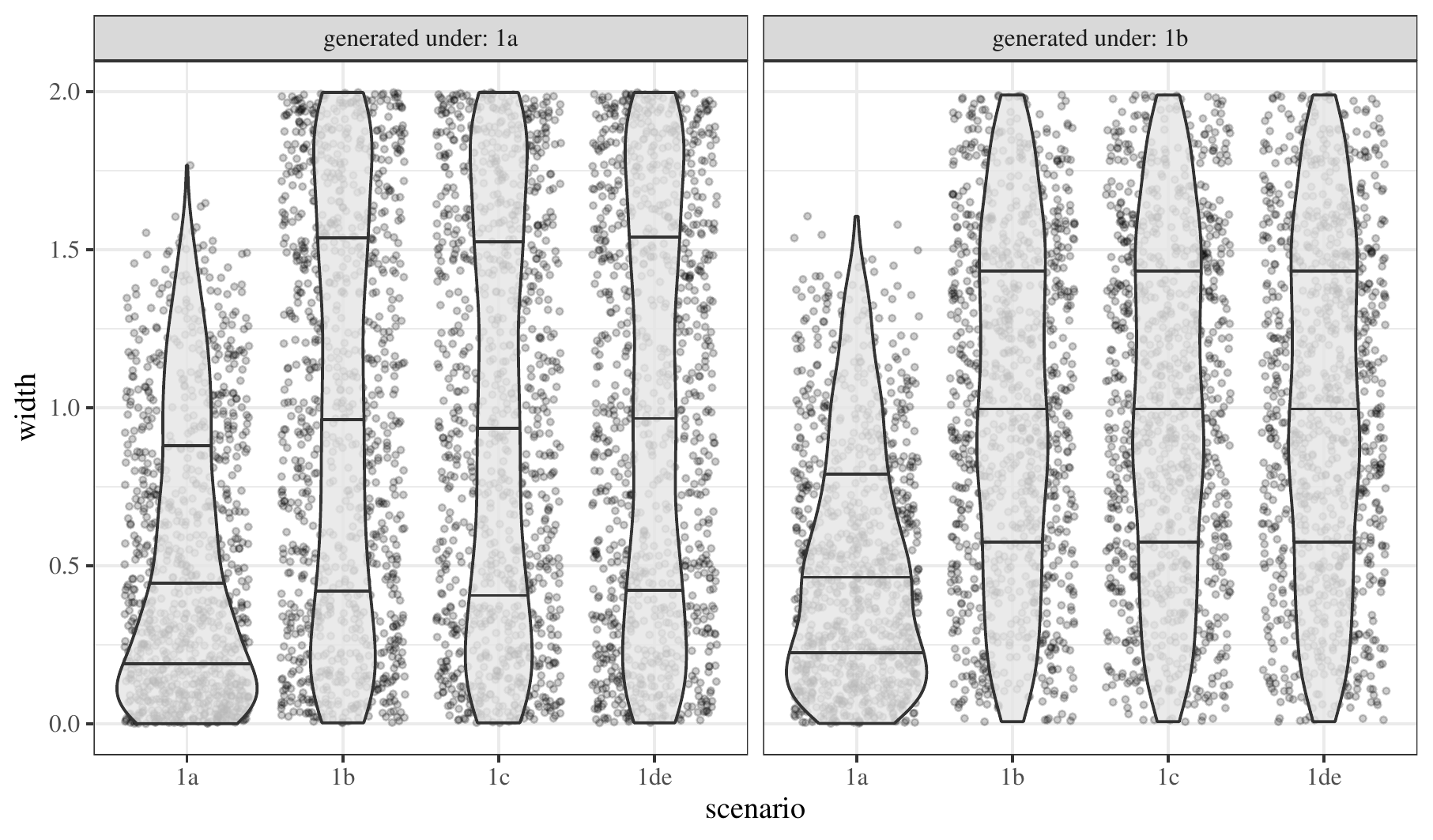}
\caption{\label{truewidth1} Comparison of the width of the true bounds for the DAGs in Figure 1 for distributions that are generated under Figure 1a (left panel) and 1b (right panel).}
\end{subfigure}

\begin{subfigure}{.95\textwidth}
\centering
\includegraphics[width = .99\textwidth]{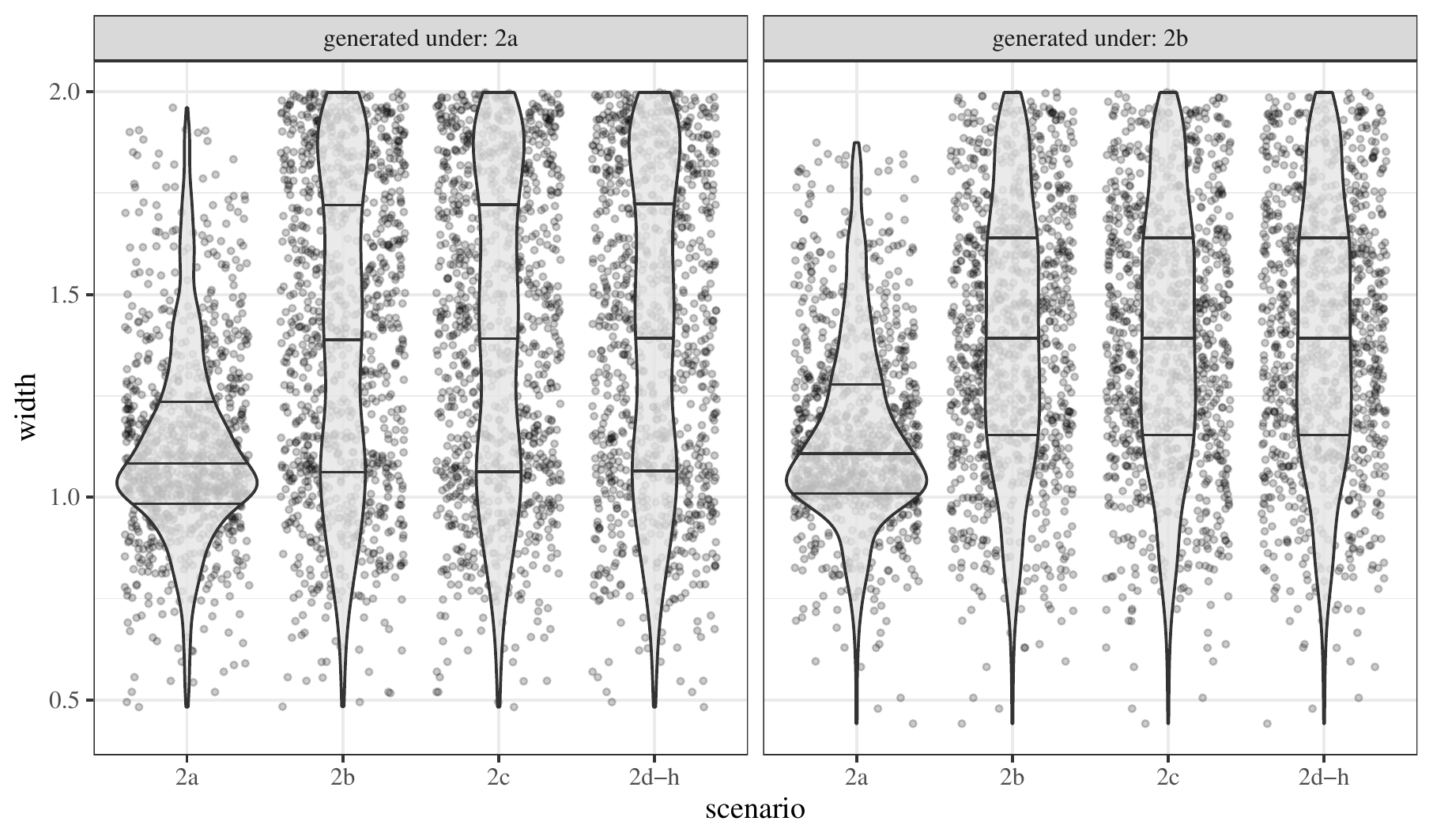}
\caption{\label{truewidth2} Comparison of the width of the true bounds for the DAGs in Figure 2 for distributions that are generated under Figure 2a (left panel) and 2b (right panel).}
\end{subfigure}
\caption{Comparison of the width of the true bounds.}
\end{figure}

\clearpage

\bibliographystyle{biometrika}
\bibliography{refer.bib}